\documentclass[12pt]{colt2018} 

\title[Solving Random Systems of Quadratic Equations with TanhWF]{Solving Random Systems of Quadratic Equations with Tanh Wirtinger
Flow}
\usepackage{times}
\usepackage{mathtools}

 \coltauthor{\Name{Zhenwei Luo} \Email{zl24@rice.edu}\and
 \Name{Ye Zhang} \Email{yz61@rice.edu}\\
 \addr {Rice University, Houston}
 }

\begin{document}

\maketitle

\begin{abstract}
    Solving quadratic systems of equations in $n$ variables and $m$
measurements of the form $y_{i}=|\boldsymbol{a}_{i}^{T}\boldsymbol{x}|^{2},i=1,...,m$
and $\boldsymbol{x}\in\mathbb{R}^{n}$, which is also known as phase
retrieval, is a hard nonconvex problem. In the case of standard Gaussian
measurement vectors, the wirtinger flow algorithm \cite{chen_solving_2015}
is an efficient solution. In this paper, we proposed a new form of
wirtinger flow and a new spectral initialization method based on this
new algorithm. We proved that the new wirtinger flow and initialization
method achieve linear sample and computational complexities. We further
extended the new phasing algorithm by combining it with other existing
methods. Finally, we demonstrated the effectiveness of our new method
in the low data to parameter ratio settings where the number of measurements
which is less than information-theoretic limit, namely, $m<2n$, via
numerical tests. For instance, our method can solve the quadratic
systems of equations with gaussian measurement vector with probability
$\ge97\%$ when $m/n=1.7$ and $n=1000$, and with probability $\approx60\%$
when $m/n=1.5$ and $n=1000$.

\end{abstract}


\section{Introduction}

Reconstructing signal from intensity measurements only, which is also
known as phase retrieval, is an important problem with applications
in various fields, including the X-ray crystallography, astronomy,
and diffraction imaging \cite{fienup_phase_1982,gerchberg_practical_1972,hauptman_phase_1991}.
There is a recent resurgence of interest in solving phase retrieval
problem in machine learning community. A typical setting of those
works is that, the observed data $y_{i}$ is of the form, 
\[
y_{i}\approx|\langle\boldsymbol{a}_{i},\boldsymbol{x}\rangle|^{2},\ i=1,...,m\ ,
\]
where $\boldsymbol{a}_{i}$ is a random gaussian measurement vector
and $\boldsymbol{x}$ is unknown \cite{candes_phaselift:_2013}. Numerous
approaches have been proposed to solve this problem. They can be mainly
categorized into two classes. The first one is converting the problem
into a convex program which can be easily solved. PhaseLift \cite{candes_phaselift:_2013}
and PhaseCut \cite{waldspurger2015phase} all fall into this category.
It has been established that this kind of algorithm can recover any
vector $\boldsymbol{x}\in\mathbb{C}^{n}$ exactly from only $O(n)$
samples \cite{candes2014solving}. However, the computational complexities
of those convex surrogates scale as $n^{4.5}\log\frac{1}{\epsilon}$
, which limits applicability to high dimensional data. Another class
of algorithm is to optimize the nonconvex problem directly. One of
the most successful algorithm in this category is the Wirtinger Flow
(WF), which was proposed by \cite{candes_phase_2015}.
Chen et al. improved WF by wisely discarding certain outlier gradients,
and named this procedure as the Truncated Wirtinger Flow (TWF) \cite{chen_solving_2015}.
The TWF method achieves linear sample complexity and linear-time computational
cost, which is the optimal statistical complexity one can obtain for
solving a well-posed problem. There are numerous follow up works focus
on improving the truncation rules \cite{kolte_phase_2016,wang_solving_2016,zhang_provable_2016,wang2017solving}.
These works lower the critical sampling ratio required to exactly solve
the random systems of quadratic equations. In realistic
scenarios, the systems to be solved often admit low data to parameter
ratio, e.g $<2$ for real-valued system. For instance, the crystals
of important large protein complexes diffract to low resolutions,
which limits the number of measurements collected in the experiment
and compromises the model accuracy \cite{poon2007normal,schroder2010super}.
Therefore, it's of great importance to develop new methods which are
capable of solving systems with low data to parameter ratio, and can
be employed in real-world problems. In this paper, we follow the second
route to solve the nonconvex phase retrieval problem directly and
design new methods which can not only improve the success rate of
solving the random systems of quadratic equations with low data to
parameter ratio but also admits optimal statistical complexity. 

Our contributions can be summarized as follows; 1) We designed a novel
data dependent nonlinear weight function to improve the regularity
of wirtinger flow in the area with large model errors. This approach
is different from the traditional truncation rules and has shown to
increase the empirical success rates in numerical simulation. 2) We
developed a new spectral initialization method which can obtain initial
solution with high correlation w.r.t the true signal. 3) By combining
our weight function with a weight function similar to the reweighted
amplitude flow (RAF) \cite{wang2017solving}, we obtained a more powerful
phasing method, RTanhWF, and achieved the highest empirical success
rates on solving the random systems of quadratic equations relative
to state-of-the-art approaches. 4) We explored a different approach to establish
 the regularity condition which is the mainstream method to demonstrate the 
 convergence of phasing method. Instead of using the net argument, we leveraged the 
 results about the supreme of empirical process and obtained tighter bounds.

The remainder of this paper is organized as follows. We first present
how our weight functions are originated and developed, and how they
can be combined with RAF in section \ref{sec:Tanh-wirtinger-Flow}.
The numerical test results of our new methods and the TWF are shown
in section \ref{sec:Numerical-Experiment}. Detailed theoretical analysis
for our methods can be found in section \ref{sec:Theoretical-Analysis-of}.
The final section is attributed to concluding remarks.

\section{Tanh wirtinger Flow\label{sec:Tanh-wirtinger-Flow}}
Throughout this paper, we will use the following notation.
We denote the true real signal as $\boldsymbol{x}\in\mathbb{R}^{n}$,
and the design matrix as $\boldsymbol{A}\equiv[\boldsymbol{a}_{1},...,\boldsymbol{a}_{m}]^{T}\in\mathbb{R}^{m\times n},$
where $\boldsymbol{a}_{i}\sim\mathcal{N}(\boldsymbol{0},\boldsymbol{I}).$
We use $O(n)$ to represent a quantity which is of the order $Cn$,
where $C$ is a constant grater than 1. $\text{sgn}$ refers to the
sign function, $\text{sgn}(x)=\frac{x}{|x|},\forall x\in\mathbb{R}$. 

We first assume that the solution and our current estimation are statistically independent.
For simplicity, we restrict our attention to the real-valued system.
Given an estimated signal $\boldsymbol{z}$, for the gaussian random measurement vector $\boldsymbol{a}_{i}\sim\mathcal{N}(\boldsymbol{0},\boldsymbol{I})$, the corresponding observation $ \boldsymbol{a}_{i}^{T}\boldsymbol{x}$ and the
estimation $\boldsymbol{a}_{i}^{T}\boldsymbol{z}$ form a bivariate guassian, which is of the form
\begin{equation}
	p(\boldsymbol{a}_{i}^{T}\boldsymbol{z}, \boldsymbol{a}_{i}^{T}\boldsymbol{x}) = \frac{1}{2\pi\sqrt{\det(\boldsymbol{K})}}
	\exp{-\frac{\boldsymbol{b}^T\boldsymbol{K}^{-1}\boldsymbol{b}}{2}},
\end{equation}
where $\boldsymbol{b} = \begin{bmatrix} \boldsymbol{x}^T\\
										\boldsymbol{z}^T
						\end{bmatrix} \boldsymbol{a}_i$, and the covariance matrix $\boldsymbol{K}$ can be written as, 
\begin{equation}
	\boldsymbol{K} = \begin{bmatrix} \|\boldsymbol{x}\|^2 & \boldsymbol{x}^T\boldsymbol{z} \\
									\boldsymbol{z}^T\boldsymbol{x} & \|\boldsymbol{z}\|^2 
					\end{bmatrix}.
\end{equation}
Denote $\boldsymbol{a}_{i}^T\boldsymbol{x}$ as $y_i$ and $\boldsymbol{a}_{i}^T\boldsymbol{z}$ as $f_i$, the joint probability can be expanded as
\begin{equation}
	p(f_i, y_i) = \frac{1}{2\pi\sqrt{\det(\boldsymbol{K})}}\exp{-\frac{\|\boldsymbol{z}\|^2 y_i^2 - 2\boldsymbol{x}^T\boldsymbol{z} f_i y_i + \|\boldsymbol{x}\|^2 f_i^2}{2\det(\boldsymbol{K})}},
\end{equation}
and the probability of observing $y_i$ is $p(y_i) = \frac{1}{\sqrt{2\pi}\|\boldsymbol{x}\|}\exp{-\frac{y_i^2}{2\|\boldsymbol{x}\|^2}}$. Thus the probability of $f_i$ conditioned on $|y_i|$ can be otained as follows,
\begin{eqnarray}
	p(f_i\mid |y_i|) & = & \frac{p(f_i, |y_i|)}{p(|y_i|)}\\
				& = & \frac{\exp{-\frac{\sigma^2\frac{y_i^2}{\|\boldsymbol{x}\|^2}+\frac{f_i^2}{\|\boldsymbol{z}\|^2}}{2(1-\sigma^2)}}}{\sqrt{2\pi(1-\sigma^2)}\|\boldsymbol{z}\|}\cosh(\frac{\sigma}{ 1-\sigma^2}\frac{f_i |y_i|}{\|\boldsymbol{z}\|\|\boldsymbol{x}\|}),
\end{eqnarray}
where $\sigma = \frac{\boldsymbol{x}^T\boldsymbol{z}}{\|\boldsymbol{x}\|\|\boldsymbol{z}\|}$, which is the correlation between the solution $\boldsymbol{z}$ and real signal $\boldsymbol{x}$.
 Given the conditional probability $p(f_i\mid |y_i|)$, we can derive a likelihood function for the estimation of $\boldsymbol{x}$, which can be obtained by maximizing the total log likelihood 
 $\sum_{i=1}^{m}\log p(f_i\mid |y_i|)$. By observing that $\boldsymbol{z}$ and $\boldsymbol{x}$ are of the unit vector forms in the likelihood function, 
we can further simplify the likelihood to obtain the following
target function, 
\begin{equation}
	\min_{\boldsymbol{z\in S^{n-1}}}\log{\sqrt{1-\sigma^2}} + \frac{1}{2m(1-\sigma^2)}\sum_{i=1}^{m}(\frac{\sigma^2y_i^2}{\|\boldsymbol{x}\|^2}+f_i^2-2(1-\sigma^2)\log\cosh(\frac{\sigma}{ 1-\sigma^2}\frac{f_i |y_i|}{\|\boldsymbol{x}\|})).\label{eq:target}
\end{equation}
It's worth noting that similar likelihood function has been derived in crystallography fields
for a long time \cite{pannu_improved_1996,murshudov_refinement_1997,bricogne1996maximum}.
General overviews for the likelihood function can be found in \cite{murshudov_refinement_1997,lunin_likelihood-based_2002}.
 The gradient of the target function \ref{eq:target} for any $\boldsymbol{z}\in S^{n-1}$
is of the form, 
\begin{equation}
	\frac{1}{2m}\sum_{i=1}^{m}\nabla l_{i}(\boldsymbol{z})=\frac{1}{m(1-\sigma^2)}\sum_{i=1}^{m}(\boldsymbol{a}_{i}^{T}\boldsymbol{z}-\frac{\sigma\sqrt{y_{i}}}{\|\boldsymbol{x}\|}\tanh\frac{\sigma\boldsymbol{a}_{i}^{T}\boldsymbol{z}\sqrt{y_{i}}}{(1-\sigma^{2})\|\boldsymbol{x}\|})\boldsymbol{a}_{i}.\label{eq:rawgrad}
\end{equation}
 We then have the corresponding gradient descent update rule,
\[
\boldsymbol{z}^{t+1}=\boldsymbol{z}^{t}-\frac{\mu}{2m}\sum_{i=1}^{m}\nabla l_{i}(\boldsymbol{z}).
\]
 The vanilla gradient descent is not the only method to update parameters.
In fact, we can incorporate our new gradient with any first order
optimization algorithm. 

To gain some empirical understandings about the new gradient, we reexamine
it in the noiseless setting, namely, $y_{i}=|\boldsymbol{a}_{i}^{T}\boldsymbol{x}|^{2}$.
Suppose $x\in S^{n-1}$, we rewrite the gradient in equation \ref{eq:rawgrad} as 
\begin{equation}
\frac{1}{2m}\sum_{i=1}^{m}\nabla l_{i}(\boldsymbol{z})=\frac{1}{m}\sum_{i=1}^{m}\boldsymbol{a}_{i}\boldsymbol{a}_{i}^{T}(\boldsymbol{z}-\sigma\boldsymbol{x}\tanh\frac{\boldsymbol{x}^{T}\boldsymbol{a}_{i}\boldsymbol{a}_{i}^{T}\boldsymbol{z}\sigma}{1-\sigma^{2}}),\label{eq:grads}
\end{equation}
where $\frac{1}{m}\sum_{i=1}^{m}\boldsymbol{a}_{i}\boldsymbol{a}_{i}^{T}\boldsymbol{x}\tanh\frac{\boldsymbol{x}^{T}\boldsymbol{a}_{i}\boldsymbol{a}_{i}^{T}\boldsymbol{z}\sigma}{1-\sigma^{2}}$
serves as an estimator for the true signal $\boldsymbol{x}$.
It's approximated as linear combinations
of measurement vectors $\boldsymbol{a}_{i}$ where each vector is
weighted by $\boldsymbol{a}_{i}^{T}\boldsymbol{x}\tanh\frac{\boldsymbol{x}^{T}\boldsymbol{a}_{i}\boldsymbol{a}_{i}^{T}\boldsymbol{z}\sigma}{1-\sigma^{2}}$. The factor $\tanh\frac{\boldsymbol{x}^{T}\boldsymbol{a}_{i}\boldsymbol{a}_{i}^{T}\boldsymbol{z}\sigma}{1-\sigma^{2}}$ is at the core of our new form of gradient since it serves to modulate the contribution of each
observation with the estimated phase. We will investigate how this factor affects estimating $\boldsymbol{x}$.
It should be noted that the correlation between true signal and estimation, $\sigma$, in the gradient \ref{eq:grads} 
remains unknown unless the true solution is revealed. It in turn requires us to design
effective methods to estimate the correlation. In the traditional
crystallographic refinement field, this factor is often estimated
by the maximum likelihood method \cite{lunin1995r}, which introduces
additional complexity. 

We encountered the problem about designing effective weight function without any knowledge of the correlation. The main principle of such weight function is to downweight
the data point where the estimated phase is wrong, thus improving the correlation between the estimated
signal $\hat{\boldsymbol{x}}$ and the true signal $\boldsymbol{x}$. 
In this paper, we 
used a geometric observation about phase retrieval problem to achieve this goal. According
to the Grothendieck's identity from the Lemma 3.6.6 in \cite{vershynin_high_2016},
the relationship between the inner products of a random gaussian vector
$\boldsymbol{a}_{i}\in\mathbb{R}^{n}$ with any fixed vectors $\boldsymbol{x},\boldsymbol{z}\in S^{n-1}$
can be understood by reducing the problem to $\mathbb{R}^{2}$. Applying the transformation, by equation \ref{eq:polar},
we have $\boldsymbol{a}_{i}^{T}\boldsymbol{x}=r\|\boldsymbol{x}\|\cos\theta,\boldsymbol{a}_{i}^{T}\boldsymbol{z}=r\|\boldsymbol{z}\|\cos(\theta-\phi)$, where $r$ is the length of projected design vector.
A few observations lead us to consider $\cos\theta\cos(\theta-\phi)=\frac{\boldsymbol{a}_{i}^{T}\boldsymbol{x}\boldsymbol{a}_{i}^{T}\boldsymbol{z}}{r^{2}}$
as a good candidate for the weight function since its average
value is smaller in the region where $\text{sgn}(\cos\theta\cos(\theta-\phi))$
is negative, namely, where the estimatied phase is incorrect. This prompts us to estimate the length of projected
design vector. As the square norm of projected design vector is given
by $r^{2}=\frac{1}{1-\cos^{2}\phi}(\frac{(\boldsymbol{a}_{i}^{T}\boldsymbol{z})^{2}}{\|\boldsymbol{z}\|^{2}}+\frac{(\boldsymbol{a}_{i}^{T}\boldsymbol{x})^{2}}{\|\boldsymbol{x}\|^{2}}-2\frac{\boldsymbol{a}_{i}^{T}\boldsymbol{x}\boldsymbol{a}_{i}^{T}\boldsymbol{z}}{\|\boldsymbol{x}\|\|\boldsymbol{z}\|}\cos\phi)$,
suppose $\boldsymbol{x}$ and $\boldsymbol{z}$ are of equal length, we can approximate it with $(|\boldsymbol{a}_{i}^{T}\boldsymbol{x}|-|\boldsymbol{a}_{i}^{T}\boldsymbol{z}|)^{2}$
and denote it as $\sigma_{i}^{2}$. The
final weight function is of the form $\tanh\frac{\boldsymbol{x}^{T}\boldsymbol{a}_{i}\boldsymbol{a}_{i}^{T}\boldsymbol{z}}{(|\boldsymbol{a}_{i}^{T}\boldsymbol{x}|-|\boldsymbol{a}_{i}^{T}\boldsymbol{z}|)^{2}}$
and the corresponding wirtinger flow is named as TanhWFQ. We further consider alternative form for the TanhWFQ. We first reorganize $\tanh\frac{\boldsymbol{x}^{T}\boldsymbol{a}_{i}\boldsymbol{a}_{i}^{T}\boldsymbol{z}}{\sigma_{i}^{2}}$
as, 
\begin{eqnarray*}
\tanh\frac{\boldsymbol{x}^{T}\boldsymbol{a}_{i}\boldsymbol{a}_{i}^{T}\boldsymbol{z}}{(|\boldsymbol{a}_{i}^{T}\boldsymbol{x}|-|\boldsymbol{a}_{i}^{T}\boldsymbol{z}|)^{2}} & = & \text{sgn}(\boldsymbol{x}^{T}\boldsymbol{a}_{i}\boldsymbol{a}_{i}^{T}\boldsymbol{z})\tanh\frac{(|\boldsymbol{a}_{i}^{T}\boldsymbol{x}|+|\boldsymbol{a}_{i}^{T}\boldsymbol{z}|)^{2}-(|\boldsymbol{a}_{i}^{T}\boldsymbol{x}|-|\boldsymbol{a}_{i}^{T}\boldsymbol{z}|)^{2}}{4(|\boldsymbol{a}_{i}^{T}\boldsymbol{x}|-|\boldsymbol{a}_{i}^{T}\boldsymbol{z}|)^{2}}\\
 & = & \text{sgn}(\boldsymbol{x}^{T}\boldsymbol{a}_{i}\boldsymbol{a}_{i}^{T}\boldsymbol{z})\tanh((\frac{|\boldsymbol{a}_{i}^{T}\boldsymbol{x}|}{|\boldsymbol{a}_{i}^{T}\boldsymbol{x}|-|\boldsymbol{a}_{i}^{T}\boldsymbol{z}|}-\frac{1}{2})^{2}-\frac{1}{4}).
\end{eqnarray*}
It is easy to identify that the weight function consists of two layers of transformations: the
first layer is a quadratic transformation about the variable $\frac{|\boldsymbol{a}_{i}^{T}\boldsymbol{x}|}{|\boldsymbol{a}_{i}^{T}\boldsymbol{x}|-|\boldsymbol{a}_{i}^{T}\boldsymbol{z}|}$,
and the second layer applies a tanh activation function to the output
of the first layer. We may replace the quadratic transformation
with an absolute function, which can be written as $\tanh(|\frac{|\boldsymbol{a}_{i}^{T}\boldsymbol{x}|}{|\boldsymbol{a}_{i}^{T}\boldsymbol{x}|-|\boldsymbol{a}_{i}^{T}\boldsymbol{z}|}-\frac{1}{2}|-\frac{1}{2})$ and changes more conservatively in certain region.
The wirtinger flow weighted by the new alternative weight function
is named as TanhWFL, while this class of wirtinger flow is called
TanhWF.

We proceed to show how the TanhWF can be further improved by combining
with the RAF \cite{wang2017solving} and further explain our choice of weight function using probability argument. The weighting scheme proposed
in RAF applies to each $\text{sign}(\boldsymbol{a}_{i}^{T}\boldsymbol{z})(|\boldsymbol{a}_{i}^{T}\boldsymbol{x}|-|\boldsymbol{a}_{i}^{T}\boldsymbol{z}|)$,
which serves as an estimator for $\boldsymbol{a}_{i}^{T}\boldsymbol{h}$,
while our weight function acts upon the estimator for $\boldsymbol{x}$
only. We postulate that our weighting scheme is complementary to the
RAF weighting scheme. We then propose a new type of wirtinger flow
by knitting together these two weighting schemes as below, 
\[
\nabla l=\frac{1}{m}\sum_{i=1}^{m}\boldsymbol{a}_{i}(|\boldsymbol{a}_{i}^{T}\boldsymbol{z}|-\sqrt{y_{i}}f(\frac{\sqrt{y_{i}}}{\sqrt{y_{i}}-|\boldsymbol{a}_{i}^{T}\boldsymbol{z}|}))g(\frac{\sqrt{y_{i}}}{\sqrt{y_{i}}-|\boldsymbol{a}_{i}^{T}\boldsymbol{z}|})\text{sgn}(\boldsymbol{a}_{i}^{T}\boldsymbol{z}),
\]
where $f(\frac{\sqrt{y_{i}}}{\sqrt{y_{i}}-|\boldsymbol{a}_{i}^{T}\boldsymbol{z}|})$
and $g(\frac{\sqrt{y_{i}}}{\sqrt{y_{i}}-|\boldsymbol{a}_{i}^{T}\boldsymbol{z}|})$
are the weight functions in TanhWF and RAF, respectively. Since our weight functions depends solely on the value of
$\frac{\sqrt{y_{i}}}{\sqrt{y_{i}}-|\boldsymbol{a}_{i}^{T}\boldsymbol{z}|}$, we should explore the connection between its value
and the credibility of the estimated phase. 
In other words, we will compare the probability of obtaining correct estimated phase and the probability of obtaining 
wrong estimated phase for a given value of statistics. We can then construct certain weight functions to downweight the data points 
with ambiguous estimated phases.
Denote $\frac{\boldsymbol{a}_{i}^{T}\boldsymbol{x}}{\boldsymbol{a}_{i}^{T}\boldsymbol{h}}$
as $u$, using the results obtained in appendix \ref{subsec:Proof-curv},
we have $u>1$ or $u<0$ when the phases of $\boldsymbol{a}_{i}^{T}\boldsymbol{z}$
and $\boldsymbol{a}_{i}^{T}\boldsymbol{x}$ are the same, and $0\le u\le1$
when the phases of $\boldsymbol{a}_{i}^{T}\boldsymbol{z}$ and $\boldsymbol{a}_{i}^{T}\boldsymbol{x}$
are different. Given that the estimated phase agrees with as the true
phase, we have $\frac{|\boldsymbol{a}_{i}^{T}\boldsymbol{x}|}{|\boldsymbol{a}_{i}^{T}\boldsymbol{x}|-|\boldsymbol{a}_{i}^{T}\boldsymbol{z}|}=u$.
By equation \ref{eq:margu}, the probability density of $u$ is of
the form 
\[
p(u)=\frac{\rho\sqrt{1-\cos^{2}\theta}}{\pi((u\rho-\cos\theta)^{2}+1-\cos^{2}\theta)},
\]
where $\rho=\frac{\|\boldsymbol{h}\|}{\|\boldsymbol{x}\|}$, which
is the relative error, and $\cos\theta=\frac{\boldsymbol{x}^{T}\boldsymbol{h}}{\|\boldsymbol{x}\|\|\boldsymbol{h}\|}$
is the correlation between error and true signal. When the phases of
$\boldsymbol{a}_{i}^{T}\boldsymbol{z}$ and $\boldsymbol{a}_{i}^{T}\boldsymbol{x}$
differs, the statistics $\frac{|\boldsymbol{a}_{i}^{T}\boldsymbol{x}|}{|\boldsymbol{a}_{i}^{T}\boldsymbol{x}|-|\boldsymbol{a}_{i}^{T}\boldsymbol{z}|}$
can be expressed as $\frac{1}{2-\frac{1}{u}},\ u\in[0,1]$. Let $x=\frac{1}{2-\frac{1}{u}}$
with $u\in[0,1]$. By change of variables, we have 
\begin{eqnarray*}
p(x) & = & \frac{\rho\sqrt{1-\cos^{2}\theta}}{\pi((\frac{x}{2x-1}\rho-\cos\theta)^{2}+1-\cos^{2}\theta)}\frac{1}{(2x-1)^{2}}\\
 & = & \frac{\rho\sqrt{1-\cos^{2}\theta}}{\pi(x^{2}(\rho^{2}-4\rho\cos\theta+4)+(2\rho\cos\theta-4)x+1)}.
\end{eqnarray*}
We then compare the probability densities of $x$ and $u$, whose
values are in the interval $(1,\infty)\cup(-\infty,0)$, on the event
that they have the same value. Suppose $p(x)>p(u)$ and $x=u=v$,
we have the inequality,
\[
v^{2}(\rho^{2}-4\rho\cos\theta+4)+(2\rho\cos\theta-4)v+1<v^{2}\rho^{2}-2\rho\cos\theta v+1,
\]
which holds for $0<v<1$ as long as $1-\rho\cos\theta>0$.
Consequently, the most ambiguous regions are near the points where $\frac{|\boldsymbol{a}_{i}^{T}\boldsymbol{x}|}{|\boldsymbol{a}_{i}^{T}\boldsymbol{x}|-|\boldsymbol{a}_{i}^{T}\boldsymbol{z}|}\approx1$
or 0. Hence, an ideal weight function
should place small weights on the gradients around these points. Our weight function in ThanWF has exactly the desired property. In the combined wirtinger flow, we make the minimum of weight function in RAF be at $u=1$, that is, 
\[
g(\frac{\sqrt{y_{i}}}{\sqrt{y_{i}}-|\boldsymbol{a}_{i}^{T}\boldsymbol{z}|})=\tanh(w_{t}|\frac{\sqrt{y_{i}}-|\boldsymbol{a}_{i}^{T}\boldsymbol{z}|}{\sqrt{y_{i}}}-1|^{2})=\tanh(w_{t}|\frac{1}{u}-1|^{2}),
\]
where $w_{t}$ is a time varying coefficient, and slightly adjust the weight function $f$ from TanhWFL. Since $g$ has already suppressed the magnitude of
gradient around the point $u=1$, we increase the value of $f$ around
this point while reducing the value of $f$ around $u=0$ by setting
it as $f(\frac{\sqrt{y_{i}}}{\sqrt{y_{i}}-|\boldsymbol{a}_{i}^{T}\boldsymbol{z}|})=\tanh(w_{t}'(|\frac{\sqrt{y_{i}}}{\sqrt{y_{i}}-|\boldsymbol{a}_{i}^{T}\boldsymbol{z}|}|+b))$,
where $w_{t}'$ is also a time varying coefficient and $b>0$ is a
bias constant. Both $w_{t}$ and $w_{t}'$ control the magnitude or the step size of
gradient. By the heuristic in \cite{candes_phase_2015}, they can
be updated using an increasing function w.r.t step, which is of the
form $w_{t}=1-w_{0}\exp(-t/T),w_{0}\in[0,1]$. Finally, we name this
method as RTanhWFL and present it with Nesterov accelerated gradient
descent \cite{nesterov_method_1983,bengio_advances_2013} in algorithm
\ref{alg:Reweighted-Tanh-Wirtinger}.

Except the update rule, the form of the tanh weighted wirtinger flow
also inspires us to propose a new initialization algorithm. Suppose
$\boldsymbol{z}=\boldsymbol{x}$, we expect $\frac{1}{m}\sum_{i=1}^{m}\boldsymbol{a}_{i}\boldsymbol{a}_{i}^{T}\boldsymbol{x}\tanh\frac{y_{i}}{\sigma^{2}}\approx\boldsymbol{x},$
namely, $\boldsymbol{x}$ is the leading eigenvector of the matrix
$\frac{1}{m}\sum_{i=1}^{m}\boldsymbol{a}_{i}\boldsymbol{a}_{i}^{T}\tanh\frac{y_{i}}{\sigma^{2}}$.
In the initialization step, we can replace the term $\sigma$
with a crude estimation. To further exclude those outliers which is
weakly correlated with $\boldsymbol{x}$, we can discard the observations
whose magnitudes are not exceeding certain threshold as it's done in \cite{wang_solving_2016,chen_solving_2015}. We summarize
the workflow of our new phasing algorithm where the initial solution is generated by the new spectral initialization algorithm, the wirtinger flow
is given by TanhWFL or TanhWFQ and the update rule is Nesterov accelerated
gradient descent \cite{nesterov_method_1983,bengio_advances_2013}
in algorithm \ref{alg:Tanh-Wirtingle-Flow}. 

\begin{algorithm}
\caption{Tanh Wirtinger Flow\label{alg:Tanh-Wirtingle-Flow}}

\SetKwInOut{Input}{Input} 
\SetKwInOut{Initialization}{Initialization} 
\SetKwInOut{Refinement}{Refinement} 
\SetKwInOut{Output}{Output}
\Input{Measurements $\{y_i|1\le i \le m\}$ and sampling vectors $\{\boldsymbol{a}_i| 1\le i \le m\}$; Initialization scale factor $\alpha$, trimming threshold $\beta$, gradient type $g$, momentum $\mu$ and step size $s$.}
\Initialization{Drawn $\boldsymbol{z}_0^0$ from $\mathcal{N}(\boldsymbol{0},\boldsymbol{I})$, and normalize it as $\boldsymbol{z}_0^0=\frac{\boldsymbol{z}_0^0}{\|\boldsymbol{z}_0^0\|}$. Set $\hat{y}=\frac{1}{m}\sum_{i=1}^m y_i$.\\ 
\For{ $t=1:T_i$}{ 
\begin{eqnarray*} \boldsymbol{z}_0^t &=& \sum_{i=1}^m \boldsymbol{a}_i\boldsymbol{a}_i^T \boldsymbol{z}_0^{t-1} \tanh \frac{y_i}{\alpha \hat{y}}\mathbb{I}(y_i > \beta \hat{y}) \\
\boldsymbol{z}_0^t &=& \frac{\boldsymbol{z}_0^t}{\|\boldsymbol{z}_0^t\|} 
\end{eqnarray*} 
}
    Set $\boldsymbol{z}_0 = \sqrt{\hat{y}}\boldsymbol{z}_0^{T_i}$.} 
\Refinement{Set $\boldsymbol{v}_0 = \boldsymbol{0}$.\\
\For{$t=1:T_r$}{  	
\begin{equation*}
w_i = \frac{\sqrt{y_i}}{\sqrt{y_i} - |\boldsymbol{a}_i^T\boldsymbol{z}_{t-1}|} - \frac{1}{2}
\end{equation*}\\
	\eIf{ $g$ = TanhWFL }{
\begin{equation*}
\nabla l_t  =  \frac{2}{m} \sum_{i=1}^m \boldsymbol{a}_i(\boldsymbol{a}_i^T\boldsymbol{z}_{t-1} - \text{sgn}(\boldsymbol{a}_i^T\boldsymbol{z}_{t-1})\sqrt{y_i}\tanh(|w_i| - \frac{1}{2}) )
\end{equation*}
}
{
\begin{equation*}
\nabla l_t  =  \frac{2}{m} \sum_{i=1}^m \boldsymbol{a}_i(\boldsymbol{a}_i^T\boldsymbol{z}_{t-1} - \text{sgn}(\boldsymbol{a}_i^T\boldsymbol{z}_{t-1})\sqrt{y_i}\tanh(w_i^2 - \frac{1}{4}) )
\end{equation*}
}
\begin{eqnarray*} 
	\boldsymbol{v}_t & = & \mu \boldsymbol{v}_{t-1} - s \nabla l_t\\ 
	\boldsymbol{z}_t & = & \boldsymbol{z}_{t-1} - \mu \boldsymbol{v}_{t-1} + (1+\mu)\boldsymbol{v}_t
\end{eqnarray*} 
}
} 
\Output{$\boldsymbol{z}_{T_r}$} 
\end{algorithm}

\begin{algorithm}

\caption{Reweighted Tanh Wirtinger Flow\label{alg:Reweighted-Tanh-Wirtinger}}

\SetKwInOut{Input}{Input} 
\SetKwInOut{Refinement}{Refinement} 
\SetKwInOut{Output}{Output}
\Input{Measurements $\{y_i|1\le i \le m\}$ and sampling vectors $\{\boldsymbol{a}_i| 1\le i \le m\}$; Initial solution $\boldsymbol{z}_0$, exponential decay parameter $T$, weights $w_f,w_g$, bias $b$, momentum $\mu$ and step size $s$.}

\Refinement{Set $\boldsymbol{v}_0 = \boldsymbol{0}$.\\
\For{$t=1:T_r$}{  	
\begin{eqnarray*}
	x_i & = & \frac{\sqrt{y_i}}{\sqrt{y_i} - |\boldsymbol{a}_i^T\boldsymbol{z}_{t-1}|}\\
	g_i & = & \tanh((1-w_ge^{-t/T})|\frac{1}{x_i} - 1|^2)\\
	f_i & = & \tanh((1-w_fe^{-t/T})(|x_i|+b))\\
\nabla l_t  &=&  \frac{2}{m} \sum_{i=1}^m g_i\boldsymbol{a}_i(\boldsymbol{a}_i^T\boldsymbol{z}_{t-1} - \text{sgn}(\boldsymbol{a}_i^T\boldsymbol{z}_{t-1})\sqrt{y_i}f_i )\\
	\boldsymbol{v}_t & = & \mu \boldsymbol{v}_{t-1} - s \nabla l_t\\ 
	\boldsymbol{z}_t & = & \boldsymbol{z}_{t-1} - \mu \boldsymbol{v}_{t-1} + (1+\mu)\boldsymbol{v}_t
\end{eqnarray*} 
}
} 
\Output{$\boldsymbol{z}_{T_r}$} 

\end{algorithm}

\section{Numerical Experiment\label{sec:Numerical-Experiment}}

\begin{figure}
\begin{centering}
\includegraphics[width=0.4\paperwidth]{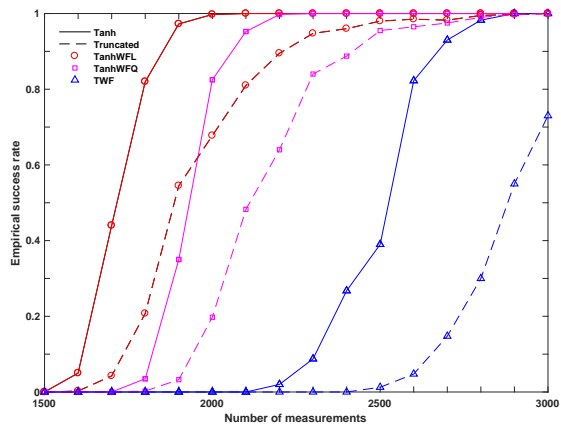}
\par\end{centering}
\caption{Empirical Success Rates on Quadratic System with $10^{3}$ Unknowns\label{fig:Empirical-Success-Rate}}
\end{figure}

In this section, we report the numerical simulation results to demonstrate
the effectiveness of our initialization method and update rules. In
all the simulations performed for TanhWF and TWF methods in this paper,
we used the following parameter settings: for the initialization stage,
we used 100 power iterations; for truncated spectral initialization,
we set the trimming threshold $\alpha_{y}=3$; for tanh weighted spectral
initialization, the scale factor $\alpha$ was set to be $4$ and
the trimming threshold $\beta$ was set to be $1$; for TWF, we adopted
the default parameters used in \cite{chen_solving_2015} when calculating
the gradient; we set the number of iterations to be $1500$, the learning
rate to be $2\times10^{-2}$ and the momentum to be $0.9$. When testing
the TWF method, we replaced the gradient calculation formula in algorithm
\ref{alg:Tanh-Wirtingle-Flow} with the formula defined in TWF. We
first compared these methods by measuring the empirical success rates
on random systems of quadratic equations with different data to parameter
ratios. The metric used to evaluate the solutions is the relative
error $\frac{\min\|\boldsymbol{z}\pm\boldsymbol{x}\|}{\|\boldsymbol{x}\|}$
which is defined in \cite{chen_solving_2015}. In these tests, we
fixed the number of unknowns to be $1000$ while varying the number
of measurements from $1500$ to $3000$ with a step size of $100$,
and we considered a system as solved if the minimum relative error
in the optimization process was smaller than $0.01$. For each method
and number of measurements, we conducted $400$ trials in which the
system to be solved was randomly generated every time and the average
values were reported. The final results for the empirical success
rate test are presented in figure \ref{fig:Empirical-Success-Rate}.
As it can be seen from figure \ref{fig:Empirical-Success-Rate}, the
tanh weighted spectral initialization significantly improved the success
rate of solving quadratic systems. Besides, the TanhWF methods also
have higher success rates comparing with TWF method when using the
same initialization method. Among all these methods, the TanhWFL with
tanh weighted spectral initialization achieves the highest empirical
success rate at every number of measurements. It can almost solve
any random system of quadratic equations (with probability $\ge99\%$)
when the sampling ratio exceeds 2. 

\begin{figure}
\begin{centering}
\includegraphics[width=0.4\paperwidth]{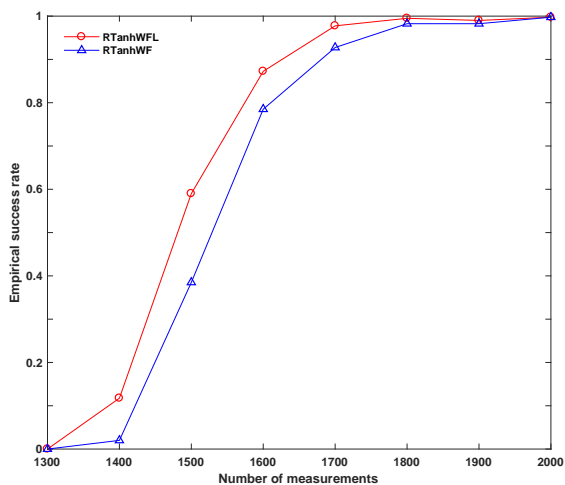}
\par\end{centering}
\caption{Empirical Success Rates for Different Reweighted Tanh Wirtinger Flows
on Quadratic System with $10^{3}$ Unknowns\label{fig:RTanhWF}}
\end{figure}

Another part of numerical experiment is to compare the empirical
success rate of RTanhWFL method and the reweighted tanh wirtinger
flow method without weighting $\sqrt{y_{i}}$, which is called RTanhWF
and implemented by letting $f=1$. We fixed the initialization method
in this test to be the tanh weighted spectral method, and didn't change
the parameter setting for this method. We still worked with the random
systems with the same number of unknowns. The parameter settings for
RTanhWFL and RTanhWF were as follows; the exponential decay parameter
$T$ was set to be $1200$; the weights $w_{f}$ and $w_{g}$ were
$0.9$ and $1$, respectively; the bias was $0.25$; the learning
rate was set to be $0.2$ and the momentum was $0.9$; the number
of iterations was 1500. All the empirical success rates reported were
the average values of 400 trials. The empirical success rates w.r.t
different numbers of measurements are shown in figure \ref{fig:RTanhWF}.
It can be seen that the RTanhWFL method can solve the random quadratic
systems with probability higher than $97\%$ when the measurement
to unknown ratio is not less than 1.7. It is also self-evident that
the RTanhWFL method outperforms the RTanhWF method at every measurement
to unknown ratio. To the best of our knowledge, this is currently
the best possible result one can obtain for solving random systems
of quadratic equations.

\begin{figure}
\begin{centering}
\includegraphics[width=0.4\paperwidth]{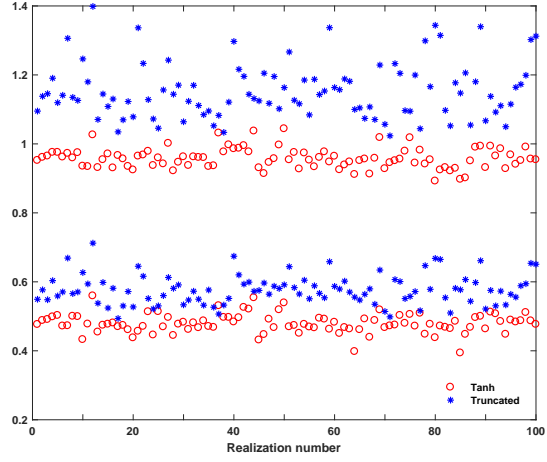}
\par\end{centering}
\caption{Relative Error and the Correlation between Error and True Solution
for Different Initialization methods\label{fig:Relative-Error-and}}
\end{figure}

We also presented the relative errors and their correlations w.r.t
true signal of the initial solutions returned by different initialization
methods. The correlation between error and true signal is defined
as 
\[
\text{corr}=\frac{\boldsymbol{x}^{T}\boldsymbol{h}}{\|\boldsymbol{x}\|\|\boldsymbol{h}\|},
\]
where $\|\boldsymbol{h}\|=\min\|\boldsymbol{z}\pm\boldsymbol{x}\|$
and $\boldsymbol{h}$ is the corresponding vector. These two quantities
for the Tanh and Truncated spectral initialization methods, which
are obtained at the sampling ratio $\frac{m}{n}=2$ and the number
of unknowns $n=1000$, are shown in figure \ref{fig:Relative-Error-and}.
The initialization errors in different realizations are shown in the
upper part of figure. It is self-evident that the initial solution
returned by the Tanh method fluctuates around 0.95, and has smaller relative error in each of
100 realizations and smaller variance. The correlation between
initial error and true signal for Tanh method fluctuates around 0.5.

\section{Convergence Analysis of the TanhWFQ method\label{sec:Theoretical-Analysis-of}}

In this section, we will perform the convergence analysis for TanhWFQ
by verifying that it satisfies the regularity condition proposed in
\cite{candes_phase_2015}. Without loss of generality, we assume $\|\boldsymbol{x}-\boldsymbol{z}\|^{2}\le\|\boldsymbol{x}+\boldsymbol{z}\|^{2}$ throughout this section.
To establish the regularity condition,
we need to bound two quantities, $-\langle\boldsymbol{h},\frac{1}{2m}\nabla l\rangle$
and $\|\frac{1}{2m}\nabla l\|$, which are the curvature and smoothness
of target function $l$, respectively. We have the following two propositions
for them.
\begin{proposition}
\label{prop:curvature}In the noiseless setting, for a fixed vector
$\boldsymbol{x}\in\mathbb{R}^{n}$, with probability at least $1-\exp[-cm\min(u,u^{2})]$,
\begin{equation}
-\langle\boldsymbol{h},\frac{1}{2m}\nabla l\rangle\ge(1-C_{3}-(C_{1}\frac{n}{m}+C_{2}\sqrt{\frac{n}{m}})-u)\|\boldsymbol{h}\|^{2},\label{eq:curvature}
\end{equation}
    holds for all $\boldsymbol{h}\in\mathbb{R}^{n}$, where $C_{1},C_{2}, C_{3}$
and $c$ are universal positive constants.
\end{proposition}

Specifically, for all $\boldsymbol{h}$ satisfying $\|\boldsymbol{h}\|\le\|\boldsymbol{x}\|$
and outside the region $\frac{\boldsymbol{x}^{T}\boldsymbol{h}}{\|\boldsymbol{x}\|\|\boldsymbol{h}\|}\in[0.6,1]\cap\frac{\|\boldsymbol{h}\|}{\|\boldsymbol{x}\|}\in[0.6,1]$,
we have $C_{3}\text{\ensuremath{\le}}0.8$. Therefore, given sufficiently
large $m/n$, it is guaranteed that $-\langle\boldsymbol{h},\frac{1}{2m}\nabla l\rangle\ge c\|\boldsymbol{h}\|^{2}$
holds for a constant $c>0$ with high probability. Moreover, the norm
of gradient satisfies the smoothness condition.
\begin{proposition}
\label{prop:smoothness}In the noiseless setting, let $\boldsymbol{A}$
be an isotropic, sub-gaussian random matrix. Given $\|\boldsymbol{A}_{i}\|_{\psi_{2}}\le K$,
there exist some universal constants $C_{1},C>0$ such that
\[
\frac{\|\nabla l\|}{\sqrt{m}}\le(C_{1}+CK^{2}[\sqrt{\frac{n}{m}}+u])\|\boldsymbol{h}\|,
\]
holds with probability at least $1-\exp[-cmu^{2}]$ for all $\boldsymbol{x},\boldsymbol{z}\in\mathbb{R}^{n}$.
\end{proposition}

The detailed proofs of proposition \ref{prop:curvature} and proposition
\ref{prop:smoothness} can be found in appendix \ref{subsec:Proof-curv}
and appendix \ref{subsec:Proof-of-smooth}. These proofs consist of two steps. In the first step, the expectations of these quantities in certain regions are calculated by assuming $\boldsymbol{z}$ is statistically independent from the measurement vectors $\boldsymbol{a}_i$. 
However, due to the fact that the expectation has no analytical solution, we used dyadic decomposition to obtain numerically integrated bounds on those regions. The second step is to demonstrate the concentration property of these
quantities for all $\boldsymbol{z}$ in those regions. Instead of using the net argument in previous works, we resort to a result about the supreme of empirical process, which is more general and provides tighter bound.

Given $-\langle\boldsymbol{h},\frac{1}{2m}\nabla l\rangle\ge c\|\boldsymbol{h}\|^{2}$
and $\|\frac{1}{2m}\nabla l\|^{2}\le C\|\boldsymbol{h}\|^{2}$ hold
with high probability, we have the following main theorem.
\begin{theorem}
\label{thm:In-the-noiseless}In the noiseless setting, for a fixed
vector $\boldsymbol{x}\in\mathbb{R}^{n}$, there exist some universal
constants $0<\rho_{0}<1$ and $c_{1},c_{2}$ such that with probability
exceeding $1-c_{1}\exp[-c_{2}m]$, 
\[
\text{dist}^{2}(\boldsymbol{z}-\frac{\mu}{m}\nabla l(\boldsymbol{z}),\boldsymbol{x})\le(1-\rho_{0})\text{dist}^{2}(\boldsymbol{z},\boldsymbol{x})
\]
 holds for all $\boldsymbol{h}$ satisfying $\|\boldsymbol{h}\|\le\|\boldsymbol{x}\|$
and not in the set $\frac{\boldsymbol{x}^{T}\boldsymbol{h}}{\|\boldsymbol{x}\|\|\boldsymbol{h}\|}\in[0.6,1]\cap\frac{\|\boldsymbol{h}\|}{\|\boldsymbol{x}\|}\in[0.6,1]$
as long as $\mu$ is sufficiently small.
\end{theorem}

Proof. Since $\|\boldsymbol{x}-\boldsymbol{z}\|^{2}\le\|\boldsymbol{x}+\boldsymbol{z}\|^{2}$, we have 
\begin{eqnarray*}
\text{dist}^{2}(\boldsymbol{z}-\frac{\mu}{m}\nabla l(\boldsymbol{z}),\boldsymbol{x}) & = & \|\boldsymbol{h}+\frac{\mu}{m}\nabla l(\boldsymbol{z})\|^{2}\\
 & = & \|\boldsymbol{h}\|^{2}+2\mu\langle\boldsymbol{h},\frac{\nabla l(\boldsymbol{z})}{m}\rangle+\mu^{2}\|\frac{\nabla l(\boldsymbol{z})}{m}\|^{2}\\
 & \le & \|\boldsymbol{h}\|^{2}(1-4c\mu+4C\mu^{2}).
\end{eqnarray*}
Therefore, as long as $1-4c\mu+4C\mu^{2}<1$, namely, $0<\mu<\frac{c}{C}$,
the gradient update is contractive and our algorithm enjoys geometric
convergence rate. The convergence properties
of TanhWFL method can also be established using similar approaches.

\subsection{Initialization via tanh weighted spectral method\label{subsec:Initialization-via-tanh}}

To demonstrate that it is possible to solve the phasing
problem with our algorithm, it remains to prove that the spectral initialization method
can return a solution which is close to the true signal $\boldsymbol{x}$.
Without loss of generality, we assume the minimum distance between
estimated solution $\boldsymbol{z}_0$ and true signal $\boldsymbol{x}$
is $\|\boldsymbol{x}-\boldsymbol{z}_{0}\|$. The distance between
$\boldsymbol{x}$ and $\boldsymbol{z}_{0}$ is bounded in the following
theorem. 
\begin{theorem}
\label{thm:spectral}Consider the model where $y_{i}=|\boldsymbol{a}_{i}^{T}\boldsymbol{x}|$
and $\boldsymbol{a}_{i}\overset{ind}{\sim}\mathcal{N}(\boldsymbol{0},\boldsymbol{I})$,
with probability exceeding $1-\exp[-cm]$, the solution $\boldsymbol{z}_{0}$
returned by the tanh weighted spectral method obeys

\[
\|\boldsymbol{z}_{0}-\boldsymbol{x}\|\le\delta\|\boldsymbol{x}\|,
\]
 where $\delta>0$ is a small constant, provided that $m>c_{0}n$
for some sufficiently large constant $c_{0}$.
\end{theorem}

Our proof starts by introducing a new vector $\boldsymbol{z}$ with
the norm $\|\boldsymbol{x}\|$ and is parallel to the vector $\boldsymbol{z}_{0}$.
Then the distance between $\boldsymbol{z}_{0}$ and $\boldsymbol{x}$
can be decomposed as,
\begin{equation}
\|\boldsymbol{z}_{0}-\boldsymbol{x}\|\le\|\boldsymbol{z}-\boldsymbol{z}_{0}\|+\|\boldsymbol{z}-\boldsymbol{x}\|.\label{eq:decomp}
\end{equation}
Without loss of generality, we can assume $\|\boldsymbol{x}\|=1$.
Since the length of $\boldsymbol{z}_{0}$ is determined by an average
of $|\boldsymbol{a}_{i}\boldsymbol{x}|^{2}$, using the inequality
154 in \cite{chen_solving_2015}, the first term in the inequality
\ref{eq:decomp} can be bounded by $\max(\sqrt{1+2\epsilon}-1,1-\sqrt{1-2\epsilon})$
with probability $1-\exp(-cm)$ for sufficiently large $m/n$. It
takes more efforts to show that the second term is also bounded by
a small constant $\epsilon'$ with high probability. The detailed
proof about the bound of second term is included in appendix.

\section{Concluding Remarks}

In this paper, we presented a new phase retrieval algorithm which
employs the tanh activation function to weight the current estimation
about the phase for each measurement. We have shown that the TanhWF
method has higher success rate in solving random systems of quadratic
equations than the TWF method when using the same initialization method
and parameter update rule. In addition, we also proposed a new tanh
weighted spectral initialization method which significantly improved
the success rate comparing with the truncated initialization method.
We have proved that the TanhWF method satisfies the regularity condition
for gaussian design matrix \cite{candes_phase_2015}. Finally, we
designed the RTanhWFL method which achieved the best possible performance
for solving random quadratic systems. It is worth pointing out that
there remain some problems to be addressed, such as completing the
convergence analysis for the RTanhWFL method and in the noisy setting,
and investigating the effect of acceleration in our method. Future
possible research extensions include extending our theoretical analysis to complex-valued
signals, and developing criteria to compare different wirtinger
flow methods. 

\acks{We thank a bunch of people.}

\bibliography{ref.bib}

\appendix
\section{Proof of proposition \ref{prop:curvature}: the local curvature condition\label{subsec:Proof-curv}}

We first verify that the curvature satisfies the lower bound, 
\[
-\langle\boldsymbol{h},\frac{1}{2m}\nabla l\rangle\ge c\|\boldsymbol{h}\|^{2},
\]
 where $c>0$ is a constant smaller than 1 in certain regions. We rewrite this quantity
to strengthen its connection with $\|\boldsymbol{h}\|$, 
\begin{eqnarray}
-\langle\boldsymbol{h},\frac{1}{2m}\nabla l\rangle & = & \frac{1}{m}\sum_{i=1}^{m}\boldsymbol{h}^{T}\boldsymbol{a}_{i}[\boldsymbol{a}_{i}^{T}\boldsymbol{x}\tanh\frac{\boldsymbol{x}^{T}\boldsymbol{a}_{i}\boldsymbol{a}_{i}^{T}\boldsymbol{z}}{\sigma_{i}^{2}}-\boldsymbol{a}_{i}^{T}\boldsymbol{z}]\nonumber \\
 & = & \frac{1}{m}\sum_{i=1}^{m}\boldsymbol{h}^{T}\boldsymbol{a}_{i}\boldsymbol{a}_{i}^{T}\boldsymbol{h}-(1-\tanh\frac{\boldsymbol{x}^{T}\boldsymbol{a}_{i}\boldsymbol{a}_{i}^{T}\boldsymbol{z}}{\sigma_{i}^{2}})\boldsymbol{h}^{T}\boldsymbol{a}_{i}\boldsymbol{a}_{i}^{T}\boldsymbol{x}.\label{eq:curv}
\end{eqnarray}
 The first term in equation \ref{eq:curv} can be bounded using the
standard result since $\boldsymbol{a}_{i}^{T}\boldsymbol{h}$ is a
simple gaussian random variable. It then boils down to showing that
$\frac{1}{m}\sum_{i=1}^{m}(1-\tanh\frac{\boldsymbol{x}^{T}\boldsymbol{a}_{i}\boldsymbol{a}_{i}^{T}\boldsymbol{z}}{\sigma_{i}^{2}})\boldsymbol{h}^{T}\boldsymbol{a}_{i}\boldsymbol{a}_{i}^{T}\boldsymbol{x}\le c\|\boldsymbol{h}\|^{2}$
holds with high probability. Our proof mainly consists of two parts:
we will calculate the expectation of this random variable and shows
it is smaller than $c\|\boldsymbol{h}\|^{2}$, where $c<1$ in certain
regions; next, we demonstrate that the sample average is concentrated
around its expectation with high probability for all $\boldsymbol{h}$
in a set. We begin with writing the random variable $(1-\tanh\frac{\boldsymbol{x}^{T}\boldsymbol{a}_{i}\boldsymbol{a}_{i}^{T}\boldsymbol{z}}{\sigma_{i}^{2}})\boldsymbol{h}^{T}\boldsymbol{a}_{i}\boldsymbol{a}_{i}^{T}\boldsymbol{x}$
as a function of $\boldsymbol{a}_{i}^{T}\boldsymbol{x}$ and $\boldsymbol{a}_{i}^{T}\boldsymbol{h}$
only. Suppose $\boldsymbol{a}_{i}^{T}\boldsymbol{x}>0$, when $\boldsymbol{a}_{i}^{T}\boldsymbol{x}$
and $\boldsymbol{a}_{i}^{T}\boldsymbol{z}$ have the same sign, we
have $\boldsymbol{a}_{i}^{T}\boldsymbol{h}<\boldsymbol{a}_{i}^{T}\boldsymbol{x}$,
which is equivalent to $\frac{\boldsymbol{a}_{i}^{T}\boldsymbol{x}}{\boldsymbol{a}_{i}^{T}\boldsymbol{h}}>1$
and $\frac{\boldsymbol{a}_{i}^{T}\boldsymbol{x}}{\boldsymbol{a}_{i}^{T}\boldsymbol{h}}<0$;
when $\boldsymbol{a}_{i}^{T}\boldsymbol{x}$ and $\boldsymbol{a}_{i}^{T}\boldsymbol{z}$
have different phases, $\boldsymbol{a}_{i}^{T}\boldsymbol{h}\ge\boldsymbol{a}_{i}^{T}\boldsymbol{x}$
always holds, thus resulting in $0\le\frac{\boldsymbol{a}_{i}^{T}\boldsymbol{x}}{\boldsymbol{a}_{i}^{T}\boldsymbol{h}}\le1$.
In the case of $\boldsymbol{a}_{i}^{T}\boldsymbol{x}\le0$, $\frac{\boldsymbol{a}_{i}^{T}\boldsymbol{x}}{\boldsymbol{a}_{i}^{T}\boldsymbol{h}}$
are also in the same ranges under these conditions. Hence, denote
$(1-\tanh\frac{\boldsymbol{x}^{T}\boldsymbol{a}_{i}\boldsymbol{a}_{i}^{T}\boldsymbol{z}}{\sigma_{i}^{2}})\boldsymbol{h}^{T}\boldsymbol{a}_{i}\boldsymbol{a}_{i}^{T}\boldsymbol{x}$
as $X_{i}$, $\boldsymbol{a}_{i}^{T}\boldsymbol{h}$ as $t$ and $\boldsymbol{a}_{i}^{T}\boldsymbol{x}$
as $s$, we can express $X_{i}$ by 
\begin{equation}
X_{i}=t^{2}f(\frac{s}{t})=t^{2}\begin{cases}
(1-\tanh((\frac{s}{t}-\frac{1}{2})^{2}-\frac{1}{4}))\frac{s}{t}, & \frac{s}{t}>1,\frac{s}{t}<0\\
(1-\tanh(\frac{1}{4}-\frac{1}{16(\frac{s}{t}-\frac{1}{2})^{2}}))\frac{s}{t}, & 0\le\frac{s}{t}\le1
\end{cases},\label{eq:casetanh}
\end{equation}
namely, the random variable $X_{i}$ is the product of a function
$f(\frac{s}{t})$ and $t^{2}$. We then turn to derive the joint distribution
of $\frac{s}{t}$ and $t$ which can be used to calculate the expectation
of $X_{i}$. We have the following proposition about the joint distribution
of $\frac{s}{t}$ and $t$.
\begin{proposition}
\label{prop:hprop}Suppose $\boldsymbol{a}_{i}\in\mathbb{R}^{n}$
is a random Gaussian vector where each element has zero mean and unit
variance, given two vectors $\boldsymbol{x},\boldsymbol{h}\in\mathbb{R}^{n}$,
the joint distribution of $\frac{\boldsymbol{a}_{i}^{T}\boldsymbol{x}}{\boldsymbol{a}_{i}^{T}\boldsymbol{h}}$
and $\boldsymbol{a}_{i}^{T}\boldsymbol{h}$ is
\begin{equation}
p(\frac{\boldsymbol{a}_{i}^{T}\boldsymbol{x}}{\boldsymbol{a}_{i}^{T}\boldsymbol{h}}=u,\boldsymbol{a}_{i}^{T}\boldsymbol{h}=t)=\frac{\exp(-(\frac{(u\|\boldsymbol{h}\|-\frac{\boldsymbol{x}^{T}\boldsymbol{h}}{\|\boldsymbol{h}\|})^{2}}{\|\boldsymbol{x}\|^{2}\|\boldsymbol{h}\|^{2}-|\boldsymbol{x}^{T}\boldsymbol{h}|^{2}}+\frac{1}{\|\boldsymbol{h}\|^{2}})\frac{t^{2}}{2})}{2\pi\sqrt{\|\boldsymbol{x}\|^{2}\|\boldsymbol{h}\|^{2}-|\boldsymbol{x}^{T}\boldsymbol{h}|^{2}}}|t|.\label{eq:propev}
\end{equation}
\end{proposition}

\subsection{Proof of the Proposition \ref{prop:hprop}}

Since the joint probability of $\boldsymbol{a}_{i}^{T}\boldsymbol{h}$
and $\boldsymbol{a}_{i}^{T}\boldsymbol{x}$ is a bivariate gaussian
distribution with covariance matrix, 
\[
\boldsymbol{\Sigma}=\begin{bmatrix}\|\boldsymbol{h}\|^{2} & \boldsymbol{x}^{T}\boldsymbol{h}\\
\boldsymbol{x}^{T}\boldsymbol{h} & \|\boldsymbol{x}\|^{2}
\end{bmatrix},
\]
 the conditional distribution of $\boldsymbol{a}_{i}^{T}\boldsymbol{x}$
given $\boldsymbol{a}_{i}^{T}\boldsymbol{h}$ \cite{jensen2000statistics}
can be expressed as, 
\[
p(\boldsymbol{a}_{i}^{T}\boldsymbol{x}=s|\boldsymbol{a}_{i}^{T}\boldsymbol{h}=t)=\frac{1}{\sqrt{2\pi}\sqrt{\|\boldsymbol{x}\|^{2}-\frac{|\boldsymbol{x}^{T}\boldsymbol{h}|^{2}}{\|\boldsymbol{h}\|^{2}}}}\exp(-\frac{(s-\frac{\boldsymbol{x}^{T}\boldsymbol{h}}{\|\boldsymbol{h}\|^{2}}t)^{2}}{2(\|\boldsymbol{x}\|^{2}-\frac{|\boldsymbol{x}^{T}\boldsymbol{h}|^{2}}{\|\boldsymbol{h}\|^{2}})}),
\]
while the marginal distribution of $\boldsymbol{a}_{i}^{T}\boldsymbol{h}$
is a gaussian with variance $\|\boldsymbol{h}\|^{2}$. Denote $\frac{\boldsymbol{a}_{i}^{T}\boldsymbol{x}}{\boldsymbol{a}_{i}^{T}\boldsymbol{h}}$
as $u$, $\boldsymbol{a}_{i}^{T}\boldsymbol{h}$ as $t$ and $\boldsymbol{a}_{i}^{T}\boldsymbol{x}$
as $s$, we then have 
\[
s=ut,\ t=t.
\]
The Jacobian matrix of this transformation is 
\[
\begin{vmatrix}\frac{\partial s}{\partial u} & \frac{\partial s}{\partial t}\\
\frac{\partial t}{\partial u} & \frac{\partial t}{\partial t}
\end{vmatrix}=\begin{vmatrix}t & u\\
0 & 1
\end{vmatrix}=|t|.
\]
Using change of variables, the joint probability of $u$ and $t$
is 
\begin{eqnarray*}
p(u,t) & = & p(s|t)p(t)|t|\\
 & = & \frac{|t|}{2\pi\sqrt{\|\boldsymbol{x}\|^{2}\|\boldsymbol{h}\|^{2}-|\boldsymbol{x}^{T}\boldsymbol{h}|^{2}}}\exp(-\frac{t^{2}(u-\frac{\boldsymbol{x}^{T}\boldsymbol{h}}{\|\boldsymbol{h}\|^{2}})^{2}}{2(\|\boldsymbol{x}\|^{2}-\frac{|\boldsymbol{x}^{T}\boldsymbol{h}|^{2}}{\|\boldsymbol{h}\|^{2}})}-\frac{t^{2}}{2\|\boldsymbol{h}\|^{2}}).
\end{eqnarray*}
We are arriving at equation \ref{eq:propev} after simplifying the
exponent. We can calculate the marginal probability of $u$ as 
\begin{eqnarray*}
p(u) & = & \int_{-\infty}^{\infty}p(u,t)dt\\
 & = & \frac{\frac{\|\boldsymbol{h}\|}{\|\boldsymbol{x}\|}\sqrt{1-\frac{|\boldsymbol{x}^{T}\boldsymbol{h}|^{2}}{\|\boldsymbol{x}\|^{2}\|\boldsymbol{h}\|^{2}}}}{\pi((u\frac{\|\boldsymbol{h}\|}{\|\boldsymbol{x}\|}-\frac{\boldsymbol{x}^{T}\boldsymbol{h}}{\|\boldsymbol{h}\|\|\boldsymbol{x}\|})^{2}+(1-\frac{|\boldsymbol{x}^{T}\boldsymbol{h}|^{2}}{\|\boldsymbol{x}\|^{2}\|\boldsymbol{h}\|^{2}}))}.
\end{eqnarray*}
Let $\rho=\frac{\|\boldsymbol{h}\|}{\|\boldsymbol{x}\|}$ and $\cos\theta=\frac{\boldsymbol{x}^{T}\boldsymbol{h}}{\|\boldsymbol{h}\|\|\boldsymbol{x}\|}$,
simplification of the above formula yields
\begin{equation}
p(u)=\frac{\rho\sqrt{1-\cos^{2}\theta}}{\pi((u\rho-\cos\theta)^{2}+1-\cos^{2}\theta)}.\label{eq:margu}
\end{equation}

\subsection{Proof of the expectation\label{subsec:Proof-of-exp}}
The expectation of $(1-\tanh\frac{\boldsymbol{x}^{T}\boldsymbol{a}_{i}\boldsymbol{a}_{i}^{T}\boldsymbol{z}}{\sigma_{i}^{2}})\boldsymbol{h}^{T}\boldsymbol{a}_{i}\boldsymbol{a}_{i}^{T}$
can be calculated as, 
\begin{eqnarray}
\mathbb{E}(X_{i}) & = & \int_{-\infty}^{\infty}\int_{-\infty}^{\infty}t^{2}f(u)p(u,t)dtdu\nonumber \\
 & = & 2\int_{-\infty}^{\infty}\int_{0}^{\infty}f(u)\frac{\exp(-(\frac{(u\|\boldsymbol{h}\|-\frac{\boldsymbol{x}^{T}\boldsymbol{h}}{\|\boldsymbol{h}\|})^{2}}{\|\boldsymbol{x}\|^{2}\|\boldsymbol{h}\|^{2}-|\boldsymbol{x}^{T}\boldsymbol{h}|^{2}}+\frac{1}{\|\boldsymbol{h}\|^{2}})\frac{t^{2}}{2})}{2\pi\sqrt{\|\boldsymbol{x}\|^{2}\|\boldsymbol{h}\|^{2}-|\boldsymbol{x}^{T}\boldsymbol{h}|^{2}}}t^{3}dtdu\nonumber \\
 & = & 2\int_{-\infty}^{\infty}f(u)\frac{(\frac{(u\|\boldsymbol{h}\|-\frac{\boldsymbol{x}^{T}\boldsymbol{h}}{\|\boldsymbol{h}\|})^{2}}{\|\boldsymbol{x}\|^{2}\|\boldsymbol{h}\|^{2}-|\boldsymbol{x}^{T}\boldsymbol{h}|^{2}}+\frac{1}{\|\boldsymbol{h}\|^{2}})^{-2}}{\pi\sqrt{\|\boldsymbol{x}\|^{2}\|\boldsymbol{h}\|^{2}-|\boldsymbol{x}^{T}\boldsymbol{h}|^{2}}}du.\label{eq:integ}
\end{eqnarray}
Let the angle between $\boldsymbol{x}$ and $\boldsymbol{h}$ be $\theta$,
and $\frac{\|\boldsymbol{h}\|}{\|\boldsymbol{x}\|}$ be $\rho$, the
factor in equation \ref{eq:integ} can be simplified as 
\begin{eqnarray*}
\frac{(\frac{(u\|\boldsymbol{h}\|-\frac{\boldsymbol{x}^{T}\boldsymbol{h}}{\|\boldsymbol{h}\|})^{2}}{\|\boldsymbol{x}\|^{2}\|\boldsymbol{h}\|^{2}-|\boldsymbol{x}^{T}\boldsymbol{h}|^{2}}+\frac{1}{\|\boldsymbol{h}\|^{2}})^{-2}}{\sqrt{\|\boldsymbol{x}\|^{2}\|\boldsymbol{h}\|^{2}-|\boldsymbol{x}^{T}\boldsymbol{h}|^{2}}} & = & \frac{(\|\boldsymbol{x}\|^{2}\|\boldsymbol{h}\|^{2}-|\boldsymbol{x}^{T}\boldsymbol{h}|^{2})^{\frac{3}{2}}}{((u\|\boldsymbol{h}\|-\frac{\boldsymbol{x}^{T}\boldsymbol{h}}{\|\boldsymbol{h}\|})^{2}+\|\boldsymbol{x}\|^{2}-\frac{|\boldsymbol{x}^{T}\boldsymbol{h}|^{2}}{\|\boldsymbol{h}\|^{2}})^{2}}\\
 & = & \frac{(1-\cos^{2}\theta)^{\frac{3}{2}}\rho\|\boldsymbol{h}\|^{2}}{((u\rho-\cos\theta)^{2}+1-\cos^{2}\theta)^{2}}.
\end{eqnarray*}
The expectation of $X_{i}$ then is
\[
\mathbb{E}(X_{i})=\int_{-\infty}^{\infty}f(u)\frac{2(1-\cos^{2}\theta)^{\frac{3}{2}}\rho\|\boldsymbol{h}\|^{2}}{\pi((u\rho-\cos\theta)^{2}+1-\cos^{2}\theta)^{2}}du.
\]
 This integral has no analytical solution. We thus continue to bound
the expectation of $X_{i}$ using dyadic decomposition and numerical
simulation. For $u>1$, we split the interval into $[u_{n-1},u_{n}]$,
where $u_{n}=\frac{1+\sqrt{2^{n+2}+1}}{2}$ and $n\in[-\infty,\infty]$.
On each dyadic interval, we have $f(u)\le(1-\tanh2^{n-1})u_{n},\ \forall u\in[u_{n-1},u_{n}]$.
Besides, we have the integral
\[
\int_{u_{n-1}}^{u_{n}}\frac{2(1-\cos^{2}\theta)^{\frac{3}{2}}\rho}{((u\rho-\cos\theta)^{2}+1-\cos^{2}\theta)^{2}}du=\frac{a_{n}}{1+a_{n}^{2}}-\frac{a_{n-1}}{1+a_{n-1}^{2}}+\arctan a_{n}-\arctan a_{n-1},
\]
where $a_{n}=\frac{1}{\sqrt{1-\cos^{2}\theta}}(-\cos\theta+\rho u_{n})$.
Hence, the expectation of $X_{i}$ on this interval, $\mathbb{E}(X_{i,n}^{+}$),
can be upper bounded by 
\begin{equation}
\frac{\mathbb{E}(X_{i,n}^{+})}{\|\boldsymbol{h}\|^{2}}\le(1-\tanh2^{n-1})\frac{u_{n}}{\pi}(\frac{a_{n}}{1+a_{n}^{2}}-\frac{a_{n-1}}{1+a_{n-1}^{2}}+\arctan a_{n}-\arctan a_{n-1}).\label{eq:bound1}
\end{equation}
For $n\ge1$, the decaying rate of the above bound is at least of
the order $e^{-2^{n}}2^{\frac{n}{2}}$ since $\frac{u}{1+u^{2}}+\arctan u$
is a bounded function and $(1-\tanh2^{n-1})u_{n}\sim O(e^{-2^{n}}2^{\frac{n}{2}})$.
If $n\le0$, $1-\tanh2^{n-1}$ will not decay exponentially. However,
since $u_{n}\le u_{0}$ for all $n\le0$, the expectation is upper
bounded by 
\begin{equation}
\frac{\mathbb{E}(X_{i,n}^{+})}{\|\boldsymbol{h}\|^{2}}\le\frac{u_{0}}{\pi}(\frac{a_{n}}{1+a_{n}^{2}}-\frac{a_{n-1}}{1+a_{n-1}^{2}}+\arctan a_{n}-\arctan a_{n-1}).\label{eq:simpb}
\end{equation}
Summing up equation \ref{eq:simpb} for all $n\le0$ gives rise to
\begin{equation}
\frac{\mathbb{E}(X^{+})}{\|\boldsymbol{h}\|^{2}}\le\frac{u_{0}}{\pi}(\frac{a_{0}}{1+a_{0}^{2}}-\frac{a_{-\infty}}{1+a_{-\infty}^{2}}+\arctan a_{0}-\arctan a_{-\infty}),\label{eq:sumb}
\end{equation}
which is a function with bounded values. For $u<0$, we construct
a set of intervals $[u_{n},u_{n-1}]$ with $u_{n}=\frac{1-\sqrt{2^{n+2}+1}}{2}$.
In this case, $f(u)\le(1-\tanh2^{n})u_{n-1}$ on the corresponding
interval. The expectation is upper bounded by 
\begin{equation}
\frac{\mathbb{E}(X_{i,n}^{-})}{\|\boldsymbol{h}\|^{2}}\le(1-\tanh2^{n})\frac{u_{n-1}}{\pi}(\frac{a_{n-1}}{1+a_{n-1}^{2}}-\frac{a_{n}}{1+a_{n}^{2}}+\arctan a_{n-1}-\arctan a_{n}).\label{eq:bound2}
\end{equation}
For $0<u<\frac{1}{2}$, we divide the interval into smaller ones $[u_{n-1},u_{n}]$
with $u_{n}=\frac{1}{2}-\frac{1}{2\sqrt{1+2^{n+2}}}$. Since $f(u)\le(1+\tanh2^{n})u_{n}$
on the corresponding interval, we have the following upper bound for
the expectation,
\[
\frac{\mathbb{E}(X_{i,n}^{-})}{\|\boldsymbol{h}\|^{2}}\le(1+\tanh2^{n})\frac{u_{n}}{\pi}(\frac{a_{n}}{1+a_{n}^{2}}-\frac{a_{n-1}}{1+a_{n-1}^{2}}+\arctan a_{n}-\arctan a_{n-1}).
\]
 For $\frac{1}{2}<u<1$, the interval is divided into $[u_{n},u_{n-1}]$
with $u_{n}=\frac{1}{2}+\frac{1}{2\sqrt{1+2^{n+2}}}$. Using the fact
that $f(u)\le(1+\tanh2^{n})\frac{u_{n-1}}{2\pi}$ results in the following
upper bound, 
\[
\frac{\mathbb{E}(X_{i,n}^{-})}{\|\boldsymbol{h}\|^{2}}\le(1+\tanh2^{n})\frac{u_{n-1}}{\pi}(\frac{a_{n-1}}{1+a_{n-1}^{2}}-\frac{a_{n}}{1+a_{n}^{2}}+\arctan a_{n-1}-\arctan a_{n}).
\]
The convergence behavior of the sum of upper bounds can be analyzed
using the same method as bounding the sum of upper bounds for $u>1$.

\begin{figure}
\begin{centering}
\includegraphics[scale=0.8]{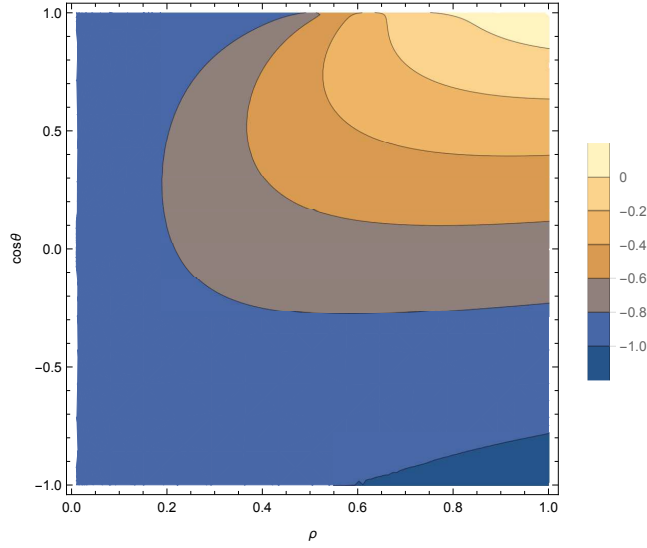}
\par\end{centering}
\caption{Contour plot for the upper bound of $\frac{\mathbb{E}(X)}{\|\boldsymbol{h}\|^{2}}-1$
w.r.t different $\frac{\|\boldsymbol{h}\|}{\|\boldsymbol{x}\|}$ and
$\cos\theta$\label{fig:Contour-plot-for}}
\end{figure}

To obtain a detailed view about how the size of bound changes with
respect to $\rho$ and $\cos\theta$, we evaluated the sum of all
upper bounds from $n=-20$ to $n=20$ over a grid with $\rho\in[0.01,1]$
and $\cos\theta\in[-0.999,0.999]$ using Mathematica. The contour
plot for the sum of upper bounds after subtracting 1 is shown in figure
\ref{fig:Contour-plot-for}. It can be seen from figure \ref{fig:Contour-plot-for}
that the upper bound of $\frac{\mathbb{E}(X)}{\|\boldsymbol{h}\|^{2}}-1$
decreases as the relative error becomes smaller and $\frac{\mathbb{E}(X)}{\|\boldsymbol{h}\|^{2}}-1<0$
in most regions. Besides, we have $\mathbb{E}((1-\tanh\frac{\boldsymbol{x}^{T}\boldsymbol{a}_{i}\boldsymbol{a}_{i}^{T}\boldsymbol{z}}{\sigma_{i}^{2}})\boldsymbol{h}^{T}\boldsymbol{a}_{i}\boldsymbol{a}_{i}^{T}\boldsymbol{x})\le0.8\|\boldsymbol{h}\|^{2}$
for all $\|\boldsymbol{h}\|,\cos\theta$ not in the set $\cos\theta\in[0.6,1]\cap\frac{\|\boldsymbol{h}\|}{\|\boldsymbol{x}\|}\in[0.6,1]$.

\subsection{Proof of the concentration property \label{subsec:sub-exp-proof}}
The next step is to demonstrate that the random variable $(1-\tanh\frac{\boldsymbol{x}^{T}\boldsymbol{a}_{i}\boldsymbol{a}_{i}^{T}\boldsymbol{z}}{\sigma_{i}^{2}})\boldsymbol{h}^{T}\boldsymbol{a}_{i}\boldsymbol{a}_{i}^{T}\boldsymbol{x}$
is tightly concentrated around its expectation. Using the formula
\ref{eq:casetanh}, $X_{i}$ is bounded by a universal constant $c$
times $t^{2}$ as, 
\begin{equation}
X_{i}\le t^{2}\underset{\frac{s}{t}}{\sup}f(\frac{s}{t})\approx1.31t^{2},\label{eq:upper-bound}
\end{equation}
 thus leading to $|(1-\tanh\frac{\boldsymbol{x}^{T}\boldsymbol{a}_{i}\boldsymbol{a}_{i}^{T}\boldsymbol{z}}{\sigma_{i}^{2}})\boldsymbol{h}^{T}\boldsymbol{a}_{i}\boldsymbol{a}_{i}^{T}\boldsymbol{x}|\le1.31(\boldsymbol{a}_{i}^{T}\boldsymbol{h})^{2}$
for all $\boldsymbol{a}_{i}^{T}\boldsymbol{h}$ and $\boldsymbol{a}_{i}^{T}\boldsymbol{x}$.
Since the sub-gaussian norm of $\boldsymbol{a}_{i}^{T}\boldsymbol{h}$
is $\sqrt{\frac{8}{3}}\|\boldsymbol{h}\|$, using lemma 2.7.6 in \cite{vershynin_high_2016}
results in the sub-exponential norm of the preceding upper bound,
namely, $\|1.31(\boldsymbol{a}_{i}^{T}\boldsymbol{h})^{2}\|_{\psi_{1}}=3.49\|\boldsymbol{h}\|^{2}$.
Moreover, we have $\mathbb{E}(\exp(|X_{i}|/C))\le\mathbb{E}(\exp(1.31(\boldsymbol{a}_{i}^{T}\boldsymbol{h})^{2}/C))$
for all $C>0$. Hence, $\|X_{i}\|_{\psi_{1}}\le1.31\|(\boldsymbol{a}_{i}\boldsymbol{h})^{2}\|_{\psi_{1}}$.
Combining all these results, we conclude that the sub-exponential
norm of $X_{i}-\mathbb{E}(X_{i})$ is smaller than $C\|\boldsymbol{h}\|^{2}$.
Using Bernstein's inequality \cite{vershynin_high_2016}, for every
$u>0$, we have 
\[
\mathbb{P}(|\frac{1}{m}\sum_{i=1}^{m}X_{i}-\mathbb{E}(X_{i})|\ge u\|\boldsymbol{h}\|^{2})\le2\exp(-mc\min(u,u^{2})),
\]
 where $c>0$ is an absolute constant, for fixed $\boldsymbol{h}$
and $\boldsymbol{x}$.

Consider $\boldsymbol{h}\in S_{\rho}^{n-1}$, which is an $n$ dimensional
Euclidean sphere with radius $\rho$, and fix $\boldsymbol{x}$, to
obtain uniform control for all $\boldsymbol{h}$ in this set $T\coloneqq S_{\rho}^{n-1}$,
it remains to show that
\[
\underset{\boldsymbol{h}\in T}{\sup}\frac{1}{m}\sum_{i=1}^{m}X_{i}-\mathbb{E}(X_{i})\ge c\|\boldsymbol{h}\|^{2}
\]
 holds with probability at most $c_{1}\exp(-c_{2}m)$ for some universal
constants $c_{1},c_{2}>0$. We can view $X_{i}$ as a random process
with the form given in \ref{eq:casetanh} and indexed by $\boldsymbol{h}$.
According to \cite{dirksen2015tail}, the supremum of the empirical
process is determined by the tail probabilities of its increments.
Hence, we begin with developing a set of properties about $X_{i}(\boldsymbol{h})$
to bound its difference. We introduce a new factor $\rho=\frac{\|\boldsymbol{h}\|}{\|\boldsymbol{x}\|}$
and normalize $\boldsymbol{h},\boldsymbol{x}$ to be unit vectors,
thus rewriting $X_{i}(\boldsymbol{h})$ as $\rho^{2}t^{2}f(\frac{s}{t\rho})\|\boldsymbol{x}\|^{2}$.
$X_{i}(\boldsymbol{h})$ now contains a scaling factor $\|\boldsymbol{x}\|^{2}$.
With a little abuse of notation, we denote $X_{i}(\boldsymbol{h})\coloneqq\frac{X_{i}(\boldsymbol{h})}{\|\boldsymbol{x}\|^{2}}$,
and turn to estimate the tail probability of 
\[
\sup_{\boldsymbol{h}\in S^{n-1}}\frac{1}{m}\sum_{i=1}^{m}X_{i}(\boldsymbol{h})-\mathbb{E}(X_{i}(\boldsymbol{h}))\ge c\rho^{2}.
\]
Let $x=\frac{s}{t\rho}$, the gradient of $\rho^{2}t^{2}f(\frac{s}{t\rho})$
w.r.t $t$ is, 
\[
\frac{\partial X_{i}(\boldsymbol{h})}{\partial t}=t\rho^{2}\begin{cases}
\frac{x^{2}}{2x-1}\text{sech}^{2}\frac{x^{2}-x}{(2x-1)^{2}}+x(1-\tanh\frac{x^{2}-x}{(2x-1)^{2}}), & 0\le x\le1\\
(2x^{3}-x^{2})\text{sech}^{2}(x^{2}-x)+x(1-\tanh(x^{2}-x)), & x<0,x>1
\end{cases}.
\]
Moreover, the second order derivate of $\rho^{2}t^{2}f(\frac{s}{t\rho})$
w.r.t $t$ is, 
\[
\frac{\partial^{2}X_{i}(\boldsymbol{h})}{\partial t^{2}}=\rho^{2}\begin{cases}
(\tanh\frac{x^{2}-x}{(2x-1)^{2}}-1)\frac{2x^{3}}{(2x-1)^{2}}\text{sech}^{2}\frac{x^{2}-x}{(2x-1)^{2}}, & 0\le x\le1\\
((2x-1)^{2}\tanh(x^{2}-x)-1)2x^{3}\text{sech}^{2}(x^{2}-x), & x<0,x>1
\end{cases}.
\]
Noting that all the functions in the case statements are bounded,
we thus conclude that $\frac{\partial X_{i}(\boldsymbol{h})}{\partial t}$
is $C\rho^{2}$-Lipstchiz w.r.t to $t$ and $|\frac{\partial X_{i}(\boldsymbol{h})}{\partial t}|\le C\rho^{2}|t|$.
Consequently, $X_{i}(\boldsymbol{h})$ is smooth; the difference between
$X_{i}(\boldsymbol{h})$ and $X_{i}(\boldsymbol{h}_{0})$ can be bounded
using Lemma 1 in \cite{sebastien_bub}, 
\begin{eqnarray*}
|X_{i}(\boldsymbol{h})-X_{i}(\boldsymbol{h}_{0})-X_{i}'(\boldsymbol{h}_{0})(\boldsymbol{a}_{i}^{T}\boldsymbol{h}-\boldsymbol{a}_{i}^{T}\boldsymbol{h}_{0})| & \le & \frac{C\rho^{2}}{2}|\boldsymbol{a}_{i}^{T}\boldsymbol{h}-\boldsymbol{a}_{i}^{T}\boldsymbol{h}_{0}|^{2},\\
|X_{i}(\boldsymbol{h})-X_{i}(\boldsymbol{h}_{0})| & \le & \rho^{2}(C_{1}|\boldsymbol{a}_{i}^{T}\boldsymbol{h}_{0}||\boldsymbol{a}_{i}^{T}(\boldsymbol{h}-\boldsymbol{h}_{0})|+C|\boldsymbol{a}_{i}^{T}(\boldsymbol{h}-\boldsymbol{h}_{0})|^{2}).
\end{eqnarray*}
We can bound the sub-exponential norm of each term in above formula
via 
\begin{eqnarray*}
\||\boldsymbol{a}_{i}^{T}\boldsymbol{h}_{0}||\boldsymbol{a}_{i}^{T}(\boldsymbol{h}-\boldsymbol{h}_{0})|\|_{\psi_{1}} & \le & \||\boldsymbol{a}_{i}^{T}\boldsymbol{h}_{0}|\|_{\psi_{2}}\||\boldsymbol{a}_{i}^{T}(\boldsymbol{h}-\boldsymbol{h}_{0})|\|_{\psi_{2}}\\
 & \le & C\|\boldsymbol{h}_{0}\|\|\boldsymbol{h}-\boldsymbol{h}_{0}\|,\\
\||\boldsymbol{a}_{i}^{T}(\boldsymbol{h}-\boldsymbol{h}_{0})|^{2}\|_{\psi_{1}} & \le & C\|\boldsymbol{h}-\boldsymbol{h}_{0}\|^{2}.
\end{eqnarray*}
 Thus, we have $\|X_{i}(\boldsymbol{h})-X_{i}(\boldsymbol{h}_{0})\|_{\psi_{1}}\le\rho^{2}(C_{1}\|\boldsymbol{h}_{0}\|+C_{2}\|\boldsymbol{h}-\boldsymbol{h}_{0}\|)\|\boldsymbol{h}-\boldsymbol{h}_{0}\|$.
With these results in hand, we can then apply the bound for the supremum
of empirical process obtained in \cite{dirksen2015tail}. Before employing
the corresponding theorem, we make a few definitions. We define a
random process $Y_{\boldsymbol{h}}$ indexed by $\boldsymbol{h}$
by setting 
\[
Y_{\boldsymbol{h}}=\frac{1}{m}\sum_{i=1}^{m}X_{i}(\boldsymbol{h})-\mathbb{E}(X_{i}(\boldsymbol{h})).
\]
Let $(T,d)$ be a metric space. A sequence $\mathcal{T}=(T_{n})_{n\ge0}$
of subsets of $T$ is called admissible if $|T_{0}|=1$ and $|T_{n}|\le2^{2^{n}}$
for all $n\ge1$. The $\gamma_{\alpha}$-functional of $(T,d)$ is
defined by 
\[
\gamma_{\alpha}(T,d)\coloneqq\underset{\mathcal{T}}{\inf}\ \underset{t\in T}{\sup}\sum_{k=0}^{\infty}2^{k/\alpha}d(t,T_{k}),\quad\forall0<\alpha<\infty,
\]
where the infimum is taken over all admissible sequences. We further
define a set of metrics as $d_{1}(s,t)=\max\|X_{t_{i}}-X_{s_{i}}\|_{\psi_{1}},$
and $d_{2}(s,t)=(\frac{1}{m}\sum_{i=1}^{m}\|X_{t_{i}}-X_{s_{i}}\|_{\psi_{1}}^{2})^{1/2}$.
\begin{lemma}
(Corollary 5.2, \cite{dirksen2015tail}) Let empirical process $Y_{t}$
be as above and let $\sigma,K>0$ be constants such that 
\[
\underset{t\in T}{\sup}\frac{1}{m}\sum_{i=1}^{m}\mathbb{E}|X_{t_{i}}-\mathbb{E}X_{t_{i}}|^{q}\le\frac{q!}{2}\sigma^{2}K^{q-2},\quad(q=2,3,...).
\]
 Then, there exist constants $c,C>0$ such that for any $u\ge1$,
\[
\mathbb{P}(\underset{t\in T}{\sup}|Y_{t}|\ge C(\frac{1}{\sqrt{m}}\gamma_{2}(T,d_{2})+\frac{1}{m}\gamma_{1}(T,d_{1}))+c(\frac{\sigma}{\sqrt{m}}\sqrt{u}+\frac{K}{m}u))\le e^{-u}.
\]
\end{lemma}

For our random process $Y_{t}$, since $\|X_{t_{i}}-\mathbb{E}X_{t_{i}}\|_{\psi_{1}}\le C\rho^{2}$
for all $t\in T$, we thus get $\sigma\le\sqrt{2}C\rho^{2}$ and $K\le C\rho^{2}$.
We next bound the $\gamma_{1}$ and $\gamma_{2}$ functionals for
$S^{n-1}$. The $\gamma_{1}$ functional of $S^{n-1}$ can be bounded
by Dudley's inequality,
\[
\gamma_{1}(T,d_{1})\le C\int_{0}^{\infty}\log N(T,d_{1},u)du,
\]
 where $N(T,d,u)$ is the covering number of $T$ of scale $u$. It
then comes down to estimating the covering number $N(T,d,u)$. Since
$d_{1}(\boldsymbol{h},\boldsymbol{h}_{0})\le\rho^{2}(C_{1}\|\boldsymbol{h}_{0}\|+C_{2}\|\boldsymbol{h}-\boldsymbol{h}_{0}\|)\|\boldsymbol{h}-\boldsymbol{h}_{0}\|$,
an $\epsilon$-net for $(S^{n-1},\|\|_{2})$ is a $\rho^{2}(C_{1}+C_{2}\epsilon)\epsilon$-net
for $(S^{n-1},d_{1})$. Namely, we have 
\[
N(S^{n-1},d_{1},\rho^{2}C\epsilon)\le N(S^{n-1},d_{1},\rho^{2}(C_{1}+C_{2}\epsilon)\epsilon)\le N(S^{n-1},\|\|_{2},\epsilon)\le(\frac{2}{\epsilon}+1)^{n}.
\]
We can thus plug the bound for covering number into Dudley's inequality,
noting that $N(T,d_{1},\rho^{2}C\epsilon)=1$ for $\epsilon\ge1$,
to get 
\begin{eqnarray*}
\gamma_{1}(T,d_{1}) & \le & C\rho^{2}\int_{0}^{1}n\log(\frac{2}{\epsilon}+1)d\epsilon\\
 & = & C_{1}\rho^{2}n.
\end{eqnarray*}

Similarly, the $\gamma_{2}$ functional can be bounded by the integral
\[
\gamma_{2}(T,d_{2})\le C\int_{0}^{\infty}\sqrt{\log N(T,d_{2},u)}du.
\]
Since $d_{2}=d_{1}$ in our case, we continue using previous bound
for the covering number. Consequently, we get 
\begin{eqnarray*}
\gamma_{2}(T,d_{2}) & \le & C\rho^{2}\int_{0}^{1}\sqrt{n\log(\frac{2}{\epsilon}+1)}d\epsilon\\
 & \le & C_{2}\rho^{2}\sqrt{n}.
\end{eqnarray*}
Combining all these results, we have the following tail bound for
the supremum of our random process, 
\begin{equation}
\mathbb{P}(\sup_{\boldsymbol{h}\in S^{n-1}}|Y_{\boldsymbol{h}}|\ge C\rho^{2}(C_{1}\frac{n}{m}+C_{2}\sqrt{\frac{n}{m}})+C_{3}\rho^{2}(\sqrt{\frac{2u}{m}}+\frac{u}{m}))\le e^{-u},\label{eq:uniform}
\end{equation}
which shows the supremum of $|Y_{\boldsymbol{h}}|$ is concentrated
within a ball of radius $C\rho^{2}$ with high probability for $m\approx n$,
thus completing the proof. 

\section{Proof of proposition \ref{prop:smoothness}: the local smoothness
condition\label{subsec:Proof-of-smooth}}

To fulfill establishing the regularity condition, we remains to verify
that
\[
\|\frac{1}{2m}\nabla l\|^{2}=\frac{1}{m^{2}}\boldsymbol{v}^{T}\boldsymbol{M}\boldsymbol{v}\le C\|\boldsymbol{h}\|^{2},
\]
 where $C$ is an absolute constant, $\boldsymbol{M}=\boldsymbol{A}^{T}\boldsymbol{A}$,
$\boldsymbol{A}\equiv[\boldsymbol{a}_{1},...,\boldsymbol{a}_{n}],$
and $\boldsymbol{v}^{T}\equiv[\boldsymbol{a}_{1}^{T}(\boldsymbol{z}-\boldsymbol{x}\tanh(\frac{\boldsymbol{x}^{T}\boldsymbol{a}_{1}\boldsymbol{a}_{1}^{T}\boldsymbol{z}}{\sigma_{i}^{2}})),...,\boldsymbol{a}_{n}^{T}(\boldsymbol{z}-\boldsymbol{x}\tanh(\frac{\boldsymbol{x}^{T}\boldsymbol{a}_{n}\boldsymbol{a}_{n}^{T}\boldsymbol{z}}{\sigma_{i}^{2}}))]$,
holds with high probability. An application of the inequality $\|\frac{1}{2m}\boldsymbol{A}^{T}\boldsymbol{v}\|\le\frac{1}{2m}\|\boldsymbol{A}\|\|\boldsymbol{v}\|$
simplifies the terms to be considered for bounding the norm of gradients.
Since we have $\|\boldsymbol{A}\|\le\sqrt{m}(1+\delta)$ from standard
result in non-asymptotic random matrix theory \cite{Vershynin2010Introduction},
the problem then boils down to controlling $\frac{1}{\sqrt{m}}\|\boldsymbol{v}\|$,
which in turn drives us to investigate the concentration property
and the expectation of $\boldsymbol{a}_{i}^{T}(\boldsymbol{z}-\boldsymbol{x}\tanh(\frac{\boldsymbol{x}^{T}\boldsymbol{a}_{i}\boldsymbol{a}_{i}^{T}\boldsymbol{z}}{\sigma_{i}^{2}}))$.
We first rewrite it to gain some insights about the connection between
$\|\boldsymbol{v}\|^{2}$ and $\|\boldsymbol{h}\|^{2}$,
\begin{eqnarray}
|\boldsymbol{a}_{i}^{T}(\boldsymbol{z}-\boldsymbol{x}\tanh(\frac{\boldsymbol{x}^{T}\boldsymbol{a}_{i}\boldsymbol{a}_{i}^{T}\boldsymbol{z}}{\sigma_{i}^{2}}))|^{2} & = & (-\boldsymbol{a}_{i}^{T}\boldsymbol{h}+\boldsymbol{a}_{i}^{T}\boldsymbol{x}(1-\tanh\frac{\boldsymbol{x}^{T}\boldsymbol{a}_{i}\boldsymbol{a}_{i}^{T}\boldsymbol{z}}{\sigma_{i}^{2}}))^{2}\nonumber \\
 & = & (\boldsymbol{a}_{i}^{T}\boldsymbol{h})^{2}-2\boldsymbol{a}_{i}^{T}\boldsymbol{h}\boldsymbol{a}_{i}^{T}\boldsymbol{x}(1-\tanh\frac{\boldsymbol{x}^{T}\boldsymbol{a}_{i}\boldsymbol{a}_{i}^{T}\boldsymbol{z}}{\sigma_{i}^{2}})\nonumber \\
 &  & +(\boldsymbol{a}_{i}^{T}\boldsymbol{x})^{2}(1-\tanh\frac{\boldsymbol{x}^{T}\boldsymbol{a}_{i}\boldsymbol{a}_{i}^{T}\boldsymbol{z}}{\sigma_{i}^{2}})^{2}.\label{eq:decomg}
\end{eqnarray}
In equation \ref{eq:decomg}, the first term $(\boldsymbol{a}_{i}^{T}\boldsymbol{h})^{2}$
is a random variable with known property, and the second term $\boldsymbol{a}_{i}^{T}\boldsymbol{h}\boldsymbol{a}_{i}^{T}\boldsymbol{x}(1-\tanh\frac{\boldsymbol{x}^{T}\boldsymbol{a}_{i}\boldsymbol{a}_{i}^{T}\boldsymbol{z}}{\sigma_{i}^{2}})$
has been investigated in previous section and has a sub-exponential
norm of size $O(\sigma^{2})$. We then proceed to show that the remaining
term $(\boldsymbol{a}_{i}^{T}\boldsymbol{x})^{2}(1-\tanh\frac{\boldsymbol{x}^{T}\boldsymbol{a}_{i}\boldsymbol{a}_{i}^{T}\boldsymbol{z}}{\sigma_{i}^{2}})^{2}$
is also a random variable with $O(\sigma^{2})$ sub-exponential norm.
We use the same method as bounding the sub-exponential norm of $\boldsymbol{h}^{T}\boldsymbol{a}_{i}\boldsymbol{a}_{i}^{T}\boldsymbol{x}(1-\tanh(\frac{\boldsymbol{x}^{T}\boldsymbol{a}_{i}\boldsymbol{a}_{i}^{T}\boldsymbol{z}}{\sigma_{i}^{2}}))$.
To show that $\mathbb{E}((\boldsymbol{a}_{i}^{T}\boldsymbol{x})^{2}(1-\tanh\frac{\boldsymbol{x}^{T}\boldsymbol{a}_{i}\boldsymbol{a}_{i}^{T}\boldsymbol{z}}{\sigma_{i}^{2}})^{2})\le C\|\boldsymbol{h}\|^{2}$,
it's enough to bound its sub-exponential norm since the expectation
of a random variable is smaller than its sub-exponential norm up to
an absolute constant \cite{vershynin_high_2016}. Of course, we can
use the method for bounding the expectation of $\boldsymbol{h}^{T}\boldsymbol{a}_{i}\boldsymbol{a}_{i}^{T}\boldsymbol{x}(1-\tanh(\frac{\boldsymbol{x}^{T}\boldsymbol{a}_{i}\boldsymbol{a}_{i}^{T}\boldsymbol{z}}{\sigma_{i}^{2}}))$
to obtain more precise bound. By equation \ref{eq:upper-bound}, we
have $|(1-\tanh\frac{\boldsymbol{x}^{T}\boldsymbol{a}_{i}\boldsymbol{a}_{i}^{T}\boldsymbol{z}}{\sigma_{i}^{2}})\boldsymbol{a}_{i}^{T}\boldsymbol{x}|\le1.31|\boldsymbol{a}_{i}^{T}\boldsymbol{h}|$
for all $\boldsymbol{a}_{i}$. Denote $(\boldsymbol{a}_{i}^{T}\boldsymbol{x})^{2}(1-\tanh\frac{\boldsymbol{x}^{T}\boldsymbol{a}_{i}\boldsymbol{a}_{i}^{T}\boldsymbol{z}}{\sigma_{i}^{2}})^{2}$
as $X_{i}$, then $X_{i}\le1.72(\boldsymbol{a}_{i}^{T}\boldsymbol{h})^{2}$
holds for all $\boldsymbol{a}_{i}$. We thus jump to the conclusion
that $X_{i}$ is a random variable with $O(\|\boldsymbol{h}\|^{2})$
sub-exponential norm and confirm that $\frac{1}{2m}\|\nabla l\|^{2}\le C\|\boldsymbol{h}\|^{2}$
holds with Bernstein type tail bounds.

The uniform bound for $\|\boldsymbol{v}\|$ in this case is easy to
obtain. For $\|\boldsymbol{A}^{T}\boldsymbol{h}\|$, according to
Theorem 1.4 in \cite{liaw2017simple}, we have, 
\begin{lemma}
(Theorem 1.4, \cite{liaw2017simple}). Let $\boldsymbol{A}$ be an
isotropic, sub-gaussian random matrix, and $T$ be a bounded subset
of $\mathbb{R}^{n}$. Given $\|\boldsymbol{A}_{i}\|_{\psi_{2}}\le K$,
for any $u\ge0$ the event 
\[
\sup_{\boldsymbol{h}\in T}|\|\boldsymbol{A}\boldsymbol{h}\|_{2}-\sqrt{m}\|\boldsymbol{h}\|_{2}|\le CK^{2}[w(T)+u\text{rad}(T)]
\]
 holds with probability at least $1-\exp(-u^{2})$. Here $\text{rad}(T)\coloneqq\sup_{\boldsymbol{h}\in T}\|\boldsymbol{h}\|_{2}$
denotes the radius of $T$. 
\end{lemma}

In our case, we set $T\coloneqq S^{n-1}$. Thus, we get $\text{rad}(T)=1$,
$w(T)\le\sqrt{n}$. The bound for the supreme of $\|\boldsymbol{A}\boldsymbol{h}\|$
can then be translated into 
\[
\sup_{\boldsymbol{h}\in S^{n-1}}\frac{1}{\sqrt{m}}|\|\boldsymbol{A}\boldsymbol{h}\|_{2}-1|\le CK^{2}[\sqrt{\frac{n}{m}}+\frac{u}{\sqrt{m}}],
\]
 which holds with probability at least $1-\exp(-u^{2})$. This immediately
implies that 
\[
\sup_{\boldsymbol{h}\in S_{\rho}^{n-1}}\frac{\|\boldsymbol{v}\|}{\sqrt{m}}\le C_{1}\rho+C\rho K^{2}[\sqrt{\frac{n}{m}}+\frac{u}{\sqrt{m}}],
\]
 holds with probability at least $1-\exp(-u^{2})$, thus completing
the proof.

\section{Proof of theorem \ref{thm:spectral}\label{subsec:Proof-of-spectral}}

This section continues bounding the second term of inequality \ref{eq:decomg}.
The vector $\boldsymbol{z}$ in the second term of inequality \ref{eq:decomg}
is the leading eigenvector of the tanh weighted design matrix, which
can be defined as 
\begin{equation}
\boldsymbol{z}=\underset{\boldsymbol{z}\in S^{n-1}}{\sup}\frac{1}{m}\sum_{i=1}^{m}\boldsymbol{z}^{T}\boldsymbol{a}_{i}\boldsymbol{a}_{i}^{T}\boldsymbol{z}\mathbb{I}(|\boldsymbol{a}_{i}^{T}\boldsymbol{x}|>\beta)\tanh\frac{|\boldsymbol{a}_{i}^{T}\boldsymbol{x}|^{2}}{\alpha},\label{eq:z_def}
\end{equation}
where $S^{n-1}$ represents the unit sphere in $\mathbb{R}^{n}$.
To show that the leading eigenvector $\boldsymbol{z}$ is close to
the true signal $\boldsymbol{x}$, we should first prove that $\boldsymbol{x}$
is the leading eigenvector of the expectation of the matrix $\frac{1}{m}\sum_{i=1}^{m}\boldsymbol{a}_{i}\boldsymbol{a}_{i}^{T}\mathbb{I}(|\boldsymbol{a}_{i}^{T}\boldsymbol{x}|>\beta)\tanh\frac{|\boldsymbol{a}_{i}^{T}\boldsymbol{x}|^{2}}{\alpha}$,
and then show that sample matrix is sufficiently close to its expectation.
With a little abuse of notation, we denote the eigenvector corresponding
to other eigenvalues as $\boldsymbol{z}$. Without loss of generality,
we assume $\boldsymbol{x}\in S^{n-1}$ since we can always absorbing
the norm of $\boldsymbol{x}$ into $\alpha,\beta$ by setting them
to be $\frac{\alpha}{\|\boldsymbol{x}\|^{2}},\frac{\beta}{\|\boldsymbol{x}\|}$.
This problem can be greatly simplified by leveraging its intrinsic
rotation invariance. $\forall\boldsymbol{x},\boldsymbol{z}\in S^{n-1}$,
we can rotate and project them with a rotation projection matrix $\boldsymbol{U}$
whose first row is $\boldsymbol{x}$, and second row is on the hyperplane
spanned by $\boldsymbol{x}$ and $\boldsymbol{z}$. Specifically,
the matrix $\boldsymbol{U}$ can be written as 
\[
\boldsymbol{U}=\begin{bmatrix}\boldsymbol{x}^{T}/\|\boldsymbol{x}\|\\
(\boldsymbol{z}-\frac{\boldsymbol{z}^{T}\boldsymbol{x}}{\|\boldsymbol{x}\|^{2}}\boldsymbol{x})^{T}/\|\boldsymbol{z}-\frac{\boldsymbol{z}^{T}\boldsymbol{x}}{\|\boldsymbol{x}\|^{2}}\boldsymbol{x}\|
\end{bmatrix}.
\]
We thus reduce the problem to $\mathbb{R}^{2}$. Applying the transform
leads to $\boldsymbol{x}'=\boldsymbol{U}\boldsymbol{x}$, whose coordinate
is $[1,0]$, and $\boldsymbol{z}'=\boldsymbol{U}\boldsymbol{z}$,
whose coordinate is $[\cos\theta,\sin\theta]$ ($\theta$ is the angle
between $\boldsymbol{x}$ and $\boldsymbol{z}$). The inner product
$\boldsymbol{a}_{i}^{T}\boldsymbol{x}$ is transformed to $\boldsymbol{a}_{i}^{T}\boldsymbol{U}^{T}\boldsymbol{x}'$,
thus prompting us to study the random vector $\boldsymbol{U}\boldsymbol{a}_{i}$.
The distribution of gaussian random vector is invariant under rotation,
while projecting the gaussian random vector from $\mathbb{R}^{n}$
to $\mathbb{R}^{2}$ yields a gaussian random vector in $\mathbb{R}^{2}.$
The 2D gaussian random vector is of the form $r[\cos\phi,\sin\phi]$,
where $\phi$ is uniformly distributed in $[0,2\pi]$, and $r$ is
distributed according to $r\exp(-\frac{r^{2}}{2})$ in $[0,\infty]$.
We can then rewrite the projections as 
\begin{equation}
\boldsymbol{a}_{i}^{T}\boldsymbol{x}=r\cos\phi,\boldsymbol{a}_{i}^{T}\boldsymbol{z}=r\cos(\phi-\theta).\label{eq:polar}
\end{equation}
 Based on the above equations, we have 
\begin{eqnarray}
(|\boldsymbol{a}_{i}^{T}\boldsymbol{x}|^{2}-|\boldsymbol{a}_{i}^{T}\boldsymbol{z}|^{2})\mathbb{I}(|\boldsymbol{a}_{i}^{T}\boldsymbol{x}|>\beta)\tanh\frac{|\boldsymbol{a}_{i}^{T}\boldsymbol{x}|^{2}}{\alpha} & = & r^{2}(\cos^{2}\phi-\cos^{2}(\phi-\theta))\nonumber \\
 &  & \mathbb{I}(|r\cos\phi|>\beta)\tanh\frac{r^{2}\cos^{2}\phi}{\alpha}.\label{eq:diff}
\end{eqnarray}
 Since the eigenvectors are orthogonal, we calculate the expectation
of \ref{eq:diff} conditioned on $\theta=\frac{\pi}{2}$. Let $f(\boldsymbol{a}_{i}^{T}\boldsymbol{x})=\mathbb{I}(|\boldsymbol{a}_{i}^{T}\boldsymbol{x}|>\beta)\tanh\frac{|\boldsymbol{a}_{i}^{T}\boldsymbol{x}|^{2}}{\alpha}$,
the expectation can be expressed as
\begin{eqnarray*}
\mathbb{E}((|\boldsymbol{a}_{i}^{T}\boldsymbol{x}|^{2}-|\boldsymbol{a}_{i}^{T}\boldsymbol{z}|^{2})f(\boldsymbol{a}_{i}^{T}\boldsymbol{x})) & = & \int_{0}^{\infty}\int_{0}^{2\pi}\frac{r^{3}}{2\pi}e^{-\frac{r^{2}}{2}}(\cos^{2}\phi-\sin^{2}\phi)\\
 &  & \mathbb{I}(|r\cos\phi|>\beta)\tanh\frac{r^{2}\cos^{2}\phi}{\alpha}d\phi dr\\
 & = & \int_{0}^{\infty}\int_{0}^{2\pi}\frac{r^{3}}{2\pi}e^{-\frac{r^{2}}{2}}(2\cos^{2}\phi-1)\\
 &  & \mathbb{I}(|r\cos\phi|>\beta)\tanh\frac{r^{2}\cos^{2}\phi}{\alpha}d\phi dr.
\end{eqnarray*}
For $2\cos^{2}\phi-1>0$, we have $\tanh\frac{r^{2}\cos^{2}\phi}{\alpha}>\tanh\frac{r^{2}}{2\alpha}$.
For $2\cos^{2}\phi-1\le0$, we have $\tanh\frac{r^{2}\cos^{2}\phi}{\alpha}\le\tanh\frac{r^{2}}{2\alpha}$.
Hence, we get $(2\cos^{2}\phi-1)\tanh\frac{r^{2}\cos^{2}}{\alpha}\ge(2\cos^{2}\phi-1)\tanh\frac{r^{2}}{2\alpha}$
for all $\phi$. Leveraging the preceding inequality leads to the
lower bound of the expectation, 
\begin{eqnarray*}
\mathbb{E}((|\boldsymbol{a}_{i}^{T}\boldsymbol{x}|^{2}-|\boldsymbol{a}_{i}^{T}\boldsymbol{z}|^{2})f(\boldsymbol{a}_{i}^{T}\boldsymbol{x})) & \ge & \int_{0}^{\infty}\int_{0}^{2\pi}\frac{r^{3}}{2\pi}e^{-\frac{r^{2}}{2}}(2\cos^{2}\phi-1)\\
 &  & \mathbb{I}(|r\cos\phi|>\beta)\tanh\frac{r^{2}}{2\alpha}d\phi dr\\
 & = & \int_{\beta}^{\infty}\frac{2r^{3}}{\pi}e^{-\frac{r^{2}}{2}}\tanh\frac{r^{2}}{2\alpha}\\
 &  & \int_{0}^{\arccos(\frac{\beta}{r})}\cos2\phi d\phi dr\\
 & = & \int_{\beta}^{\infty}\frac{2\beta r}{\pi}\sqrt{r^{2}-\beta^{2}}\tanh\frac{r^{2}}{2\alpha}e^{-\frac{r^{2}}{2}}dr\\
 & \ge & \int_{\beta}^{\infty}\frac{\beta r}{\pi}\sqrt{r^{2}-\beta^{2}}(1-e^{-\frac{r^{2}}{\alpha}})e^{-\frac{r^{2}}{2}}dr\\
 & = & \frac{\beta e^{-\frac{\beta^{2}}{2}}}{\sqrt{2\pi}}(1-e^{-\frac{\beta^{2}}{\alpha}}(\frac{\alpha}{2+\alpha})^{\frac{3}{2}})
\end{eqnarray*}
It's then easy to see that the expectation of $(|\boldsymbol{a}_{i}^{T}\boldsymbol{x}|^{2}-|\boldsymbol{a}_{i}^{T}\boldsymbol{z}|^{2})f(\boldsymbol{a}_{i}^{T}\boldsymbol{x})$
is positive. We thus confirm that $\boldsymbol{x}$ is the leading
eigenvector for $\mathbb{E}(\boldsymbol{a}_{i}\boldsymbol{a}_{i}^{T}f(\frac{|\boldsymbol{a}_{i}^{T}\boldsymbol{x}|^{2}}{\alpha}))$.
In addition, the eigengap between the largest eigenvalue and other
eigenvalues is greater than $\frac{\beta e^{-\frac{\beta^{2}}{2}}}{\sqrt{2\pi}}(1-e^{-\frac{\beta^{2}}{\alpha}}(\frac{\alpha}{2+\alpha})^{\frac{3}{2}})$.
It remains to bound the largest eigenvalue of $\mathbb{E}(\boldsymbol{a}_{i}\boldsymbol{a}_{i}^{T}f(\boldsymbol{a}_{i}^{T}\boldsymbol{x}))$,
that is, $\mathbb{E}(|\boldsymbol{a}_{i}^{T}\boldsymbol{x}|^{2}f(\boldsymbol{a}_{i}^{T}\boldsymbol{x}))$,
which can be bounded by
\begin{eqnarray*}
\mathbb{E}(|\boldsymbol{a}_{i}^{T}\boldsymbol{x}|^{2}f(\boldsymbol{a}_{i}^{T}\boldsymbol{x})) & = & \int_{0}^{\infty}\int_{0}^{2\pi}\frac{r^{3}}{2\pi}e^{-\frac{r^{2}}{2}}\mathbb{I}(|r\cos\phi|>\beta)\\
 &  & \cos^{2}\phi\tanh\frac{r^{2}\cos^{2}\phi}{\alpha}d\phi dr\\
 & \le & \int_{0}^{\infty}\int_{0}^{2\pi}\frac{r^{3}}{2\pi}e^{-\frac{r^{2}}{2}}\cos^{2}\phi\tanh\frac{r^{2}\cos^{2}\phi}{\alpha}d\phi dr\\
 & \le & \int_{0}^{\infty}\int_{0}^{2\pi}\frac{r^{3}}{2\pi}e^{-\frac{r^{2}}{2}}\cos^{2}\phi(1-\exp-\frac{2r^{2}\cos^{2}\phi}{\alpha})d\phi dr\\
 & = & \int_{0}^{\infty}\frac{r^{3}}{2}e^{-\frac{r^{2}}{2}}(1+e^{-\frac{r^{2}}{\alpha}}(I_{1}(\frac{r^{2}}{\alpha})-I_{0}(\frac{r^{2}}{\alpha})))dr\\
 & = & 1-(\frac{\alpha}{4+\alpha})^{\frac{3}{2}},
\end{eqnarray*}
where $I_{1}$ and $I_{0}$ are the modified Bessel functions of order
one and order zero, respectively. Since $f(\boldsymbol{a}_{i}^{T}\boldsymbol{x})$
is a bounded function, the subgaussian norm of the random vector $\boldsymbol{a}_{i}\sqrt{f(\boldsymbol{a}_{i}^{T}\boldsymbol{x})}$
is $O(1)$. Applying standard results on random matrices with non-isotropic
sub-gaussian rows results in the following inequality
\[
\|\frac{1}{m}\sum_{i=1}^{m}\boldsymbol{a}_{i}\boldsymbol{a}_{i}^{T}f(\frac{|\boldsymbol{a}_{i}^{T}\boldsymbol{x}|^{2}}{\alpha})-\mathbb{E}(\boldsymbol{a}_{i}\boldsymbol{a}_{i}^{T}f(\frac{|\boldsymbol{a}_{i}^{T}\boldsymbol{x}|^{2}}{\alpha}))\|\le C(1-(\frac{\alpha}{4+\alpha})^{3/2})(\sqrt{\frac{n+u}{m}}+\frac{n+u}{m}),
\]
 which holds with probability at least $1-2e^{-u}$ \cite{vershynin_high_2016}.
We denote the vector $\boldsymbol{z}$ as the solution which is defined
in equation \ref{eq:z_def}. We can then use the Davis-Kahan theorem
to deduce that 
\[
\sin\Theta(\boldsymbol{x},\boldsymbol{z})\le C\frac{\lambda_{1}}{\lambda_{1}-\lambda_{2}}(\sqrt{\frac{n+u}{m}}+\frac{n+u}{m})
\]
 holds with probability at least $1-2e^{-u}$, where $\Theta(\boldsymbol{x},\boldsymbol{z})$
is the angle between $\boldsymbol{x}$ and $\boldsymbol{z}$ \cite{yu2014useful},
$\lambda_{1},\lambda_{2}$ are the largest and second largest eigenvalues
of the expected matrix, and $\frac{\lambda_{1}}{\lambda_{1}-\lambda_{2}}\le\frac{1-(\frac{\alpha}{4+\alpha})^{\frac{3}{2}}}{\beta(1-e^{-\frac{\beta^{2}}{\alpha}}(\frac{\alpha}{2+\alpha})^{\frac{3}{2}})}\sqrt{2\pi}e^{\frac{\beta^{2}}{2}}$.
Using inequality (154) in \cite{chen_solving_2015}, we can obtain
similar bound for $\text{dist}(\boldsymbol{z}_{0},\boldsymbol{x})$.
Therefore, the estimation $\boldsymbol{z}_{0}$ returned by our spectral
initialization method will be sufficiently close to the true signal
$\boldsymbol{x}$ with probability $1-\exp[-cm]$ for $m\ge Cn$,
where $C$ is a sufficiently large constant.

\end{document}